\documentclass[11pt]{article}

\usepackage[margin=1.3in]{geometry}
\usepackage{mathtools}
\usepackage[dvips]{graphics}
\usepackage{epsfig}
\usepackage{multirow}
\usepackage{array}
\usepackage{extarrows}
\usepackage{graphicx}
\mathtoolsset{showonlyrefs}

\usepackage[symbol]{footmisc}

\usepackage{comment}
\usepackage[titletoc,title]{appendix}
\numberwithin{equation}{section}
\usepackage{amsthm}
\usepackage{enumerate}
\usepackage{amsmath}
\usepackage{amssymb,amsmath}
\usepackage{latexsym}
\usepackage[usenames]{xcolor}
\usepackage{pgf}
\usepackage[numbers]{natbib}
\usepackage{hyperref}
\hypersetup{
    colorlinks,%
    citecolor=black,%
    filecolor=black,%
    linkcolor=black,%
    urlcolor=black
}
\usepackage{subfig}
\numberwithin{equation}{section}

\newtheorem{theorem}{ \noindent T{\footnotesize HEOREM}}

\newtheorem{lemma}{ \noindent L{\footnotesize EMMA}}[section]
\newtheorem{coro}{ \noindent C{\footnotesize OROLLARY}}
\newcommand{\kw}[1]{\textcolor{magenta}{[KW: #1]}}

\newcommand{\E}{\mathbb{E}}

\newcommand{\eps}{\varepsilon}

\begin{document}

\title{Sparsity-adaptive concentration inequalities for random polynomials}
\author{Guozheng Dai$^*$ and Ke Wang$^\dagger$}

\date{\today}
\maketitle

\begin{abstract}
We prove concentration inequalities for polynomials of independent, sparse
$\alpha$-sub-exponential random variables. Specifically, we consider
$X_i=\delta_i\xi_i$, where the Bernoulli selectors $\delta_i$ are independent
with parameters $p_i$, and the variables $\xi_i$ are independent
\(\alpha\)-sub-exponential random variables (not necessarily centered). For any
polynomial $f:\mathbb R^n\to\mathbb R $ of degree at most $D$ and any
$0<\alpha \le 1 $, we establish an $L_r$-moment bound for
\(f(X)-\mathbb E f(X)\) in terms of partition norms of sparsity-weighted expected
derivative tensors. The weights count distinct coordinates rather than
multiplicities and therefore distinguish diagonal, partially diagonal, and
off-diagonal contributions. This captures the sparse scaling in both
collective fluctuation regimes and extreme-coordinate regimes.

When all sparsity parameters are equal to one, our result recovers the
polynomial concentration inequality of Götze, Sambale, and Sinulis. In degree
two, it recovers sparse Hanson-Wright bounds. As applications, we derive
deviation inequalities for the distance between a sparse simple random tensor
and a fixed subspace, and obtain lower bounds for the smallest singular value
of matrices whose columns are independent sparse simple random tensors.
\end{abstract}

\footnotetext[1]{Department of Mathematics, Hong Kong University of Science and Technology, Clear Water Bay, Kowloon, Hong Kong, guozhengdai@ust.hk. }
\footnotetext[2]{Department of Mathematics, Hong Kong University of Science and Technology, Clear Water Bay, Kowloon, Hong Kong, kewang@ust.hk. }

\bigskip

\section{Introduction}\label{sec:intro}
Concentration inequalities are a fundamental tool in probability theory and its applications. A classical example is the Gaussian concentration inequality \cite{Borell_invention}: if $G$ is a standard Gaussian vector in $\mathbb{R}^n$ and $f:\mathbb{R}^n\to \mathbb{R}$ is $1$-Lipschitz, then for every $t>0$,
\[
\mathbb{P}\bigl\{|f(G)-\mathbb{E}f(G)|\ge t\bigr\}\le 2e^{-t^2/2}.
\]
Such inequalities play a central role in the study of Gaussian processes, high-dimensional
geometry, random matrices, and related areas. However, many natural functions arising in probability, combinatorics, and theoretical computer science are not Lipschitz. A basic class of such examples is given by polynomials of independent random variables. For instance, subgraph counts in random graphs can be written as polynomial functions of independent edge indicators;
their concentration and large-deviation behavior have been extensively studied
(see, for instance, \cite{kim2000concentration,chatterjee2012missing,chatterjee2016introduction,LZ17}). 

This paper studies concentration inequalities for random polynomials with sparse independent inputs. Since tail bounds can often be obtained from $L_r$-moment estimates through Markov's inequality, much of the literature on polynomial concentration is formulated in terms of sharp bounds for $L_r$-norms; see, for example, Proposition~3.3 in \cite{gotze2021concentration}. We first recall several lines of work most relevant to the present paper.

The most classical case is the linear form \(\sum_{i=1}^{n}a_{i}\xi_{i}\), where \(a=(a_{i})_{i\le n}\) is a fixed vector and \(\{\xi_{i}\}_{i=1}^{n}\) is a sequence of independent random variables.  Montgomery-Smith \cite{MS_PAMS} gave two-sided concentration bounds for \(\sum_{i=1}^{n}a_{i}\xi_{i}\) when $\xi_i$ are Rademacher variables (i.e. $\xi_i=\pm 1$ with equal probability), Gluskin and Kwapie\'{n} \cite{gluskin1995tail} obtained the optimal \(L_{r}\)-norm estimate under log-concave tails, and Hitczenko, Montgomery-Smith and Oleszkiewicz \cite{HMO97} derived the optimal bound for the log-convex-tailed case (i.e. $\log \mathbb{P}\{|\xi_i|\ge t\}$ is convex or concave, respectively).  A general optimal estimate for
$ \|\sum_i a_i\xi_i\|_{L_r}$, under the assumption that the variables are symmetric or non-negative, was obtained by Lata{\l}a \cite{Latala97}. These results provide a mature understanding of the linear regime and serve as a benchmark for higher-order inequalities.

The next fundamental case is the quadratic form $\xi^\top A\xi-\mathbb E \xi^\top A\xi,$
which is the subject of Hanson-Wright type inequalities. The classical Hanson-Wright inequality \cite{HW71,RV13} gives a two-level tail bound governed by the Hilbert--Schmidt and operator norms of $A$. It has become a standard tool in high-dimensional probability, with applications to random
matrices, covariance estimation, compressed sensing, randomized linear algebra, graph theory, and
statistical learning. Considerable work has been devoted to extensions beyond the classical Rademacher or sub-Gaussian settings, including quadratic forms of dependent variables, non-isotropic random vectors, unbounded entries, and tensor-valued analogues; see, for example, \cite{latala99, VW15, Adamczak15,  KZ20, LC21, CY21, kim2000concentration,  gotze2021concentration,buterus2023some} and references therein. 

Sparse quadratic forms with entries $X_i=\delta_i\xi_i$ where $\delta_i\sim {\rm Bernoulli}(p_i)$ have attracted particular attention recently \cite{GLZ00,Zhou19,Zhou20,HKM19,schudy2012concentration,schudy2011bernstein,He_arXiv_HW,DHWZ_2025} (we refer to \cite{He_arXiv_HW} for a detailed discussion), as dense Hanson-Wright inequalities lose the correct dependence on sparsity parameters $p_i$. The authors of \cite{He_arXiv_HW} developed Hanson-Wright inequalities for sparse independent $\alpha$-sub-exponential variables with applications to random matrices and norm concentration. Subsequently, \cite{DHWZ_2025} obtained a refined sparse Hanson-Wright inequality for $0< \alpha \le 1$. The present paper extends this improved sparse quadratic theory to arbitrary fixed degree.

Beyond quadratic forms, a natural object is the order-$d$ tetrahedral chaos
\[\sum_{i_1,\dots,i_d} a_{i_1\dots i_d}\xi_{i_1}\dots\xi_{i_d},\]
where the $d$-tensor $A=(a_{i_1\dots i_d})_{1\le i_1,\dots,i_d\le n}$ is diagonal free (i.e., $a_{i_1\dots i_d}=0$ whenever two indices coincide) and $\xi=\{\xi_i, {i\le n}\}$ is a sequence of independent mean-zero random variables. Optimal $L_r$-norms for this chaos have been established by Lata{\l}a \cite{Latala2006} in the Gaussian case, by Adamczak and Lata{\l}a \cite{AL12} under log-concave tails (sharp for $d\le 3$), and by Kolesko and Lata{\l}a \cite{KL15} under log-convex tails. %Crucially, this literature universally assumes centered random variables. While centering by subtracting the mean is conceptually straightforward, converting the resulting bounds into neat, ready-to-use formulas is surprisingly complex, with computational difficulty growing steeply in the degree $d$. This challenge is already apparent for $d=3$ and becomes prohibitive for higher orders. \kw{Add more references here: Adamczak and Latala considered log-concave distribution; Vershynin considered subgaussian distribution with optimal dependence on the degree of the tensors, of which Bamberger, Krahmer, Ward  gave improved dependence on the dimension of the component vectors; and Adamczak, Latala, and Melle extended the result to values in some Banach spaces. St´ephane Boucheron et al. “Moment inequalities for functions of independent random variables”}

Several related developments address tensor-valued or structured higher-order models. Vershynin~\cite{vershynin2020concentration} proved concentration inequalities for simple random tensors with independent sub-Gaussian entries, with explicit dependence on both dimension and tensor order. Bamberger, Krahmer and Ward~\cite{BKW22} established a Hanson-Wright inequality for random tensors. Adamczak, Lata{\l}a and Meller~\cite{ALM21} obtained moment and tail estimates for Gaussian chaoses with values in Banach spaces. The general moment inequalities of Boucheron, Bousquet, Lugosi and Massart~\cite{BBLM05} apply to broad classes of functions of independent random variables and include moment inequalities for Rademacher chaoses.

These results primarily address either centered chaoses, dense tensor models, or general functions through abstract difference operators. While centering by subtracting the mean is conceptually straightforward, converting the resulting bounds into neat formulas is surprisingly complex, with computational difficulty growing steeply in the degree $d$. This challenge is already apparent for $d=3$ and becomes prohibitive for higher orders.

A broader framework for polynomial concentration was developed by Adamczak-Wolff \cite{AW15} and by G\"otze-Sambale-Sinulis \cite{gotze2021concentration}. These works study general polynomials $f(\xi_1,\ldots,\xi_n)$ of independent, not necessarily centered, random variables, with bounds expressed in terms of partition norms of the expected derivative tensors $\mathbb E f^{(d)}(\xi)$ for $1\le d\le D$. Adamczak and Wolff \cite{AW15} treated the sub-Gaussian setting, while G\"otze, Sambale and Sinulis \cite{gotze2021concentration} obtained corresponding inequalities for $\alpha$-sub-exponential variables with $0<\alpha\le1$. The appearance of all derivative orders is essential: lower-order derivatives may dominate fluctuations even for high-degree polynomials.

The present paper develops a sparse analogue of this polynomial concentration theory. We consider independent random variables $X_i=\delta_i\xi_i$ for $i=1,\ldots,n$, where the Bernoulli selectors $\delta_i$ have parameters $p_i$ and the $\xi_i$ are independent $\alpha$-sub-exponential with $0<\alpha\le1$. For a polynomial $f:\mathbb R^n\to\mathbb R$ of degree at most $D$, we prove moment and tail bounds for $f(X)-\mathbb E f(X)$ in terms of sparsity-weighted versions of the expected derivative tensors. The case $\alpha \ge 1$ requires different analytical tools and will be addressed in future work.

We also apply our results to sparse simple random tensors, obtaining deviation inequalities for the distance between a sparse simple random tensor and a fixed subspace. Using leave-one-out arguments and the negative second moment identity, this yields lower bounds for the smallest singular value of matrices with independent sparse simple random tensor columns.

\medskip

The rest of the paper is organized as follows. In Section~\ref{sec:main}, we state the main sparse polynomial concentration theorem and its tail consequence. Section~\ref{sec:applications} presents applications to sparse simple random tensors, including distance estimates to fixed subspaces and lower bounds for the smallest singular value of random tensor matrices. Section~\ref{sec:proofs} is devoted to the proofs. The appendix contains auxiliary tensor-norm estimates needed for the application section.

\section{Main results}\label{sec:main}

We study concentration inequalities for polynomials of independent sparse
$\alpha$-sub-exponential random variables. Let $X_i=\delta_i\xi_i$ for  $1\le i \le n,$ 
where the variables $\{\delta_i,\xi_i: i \in [n] \}$ are mutually independent,
$\delta_i \sim{\rm Bernoulli}(p_i)$ with $0< p_i \le 1$, and $\xi_i$ is
$\alpha$-sub-exponential. Our goal is to bound
\[
f(X)-\mathbb E f(X),
\qquad X=(X_1,\ldots,X_n),
\]
where $f: \mathbb R^n\to\mathbb R$ is a polynomial of degree $ D$.

Recall that a random variable $ \xi $ is called \emph{$ \alpha$-sub-exponential} if, for all $t\ge 0$,
\begin{align}\label{Eq_definition_alphasub}
\mathbb P\{|\xi|\ge t\}\le C e^{-c t^\alpha},
\end{align}
where $C,c>0$. We use the standard Orlicz functional
\[
\|\xi\|_{\Psi_\alpha} := \inf\left\{ t>0: \mathbb E\exp\left(\frac{|\xi|^\alpha}{t^\alpha}\right)\le 2 \right\}.
\]
For $\alpha \ge1$, this is a norm; for $0<\alpha<1 $, it is no longer a norm, but we keep the same notation.

Since $ |X_i|\le |\xi_i| $, one has
\[
\|X_i\|_{\Psi_\alpha}\le \|\xi_i\|_{\Psi_\alpha}.
\]
Thus the dense polynomial concentration inequalities of Adamczak-Wolff \cite{AW15} and
G\" otze-Sambale-Sinulis \cite{gotze2021concentration} may be applied directly to $f(X)-\mathbb E f(X)$. 
However, this direct approach ignores the fact that $X_i$ vanishes with probability
$1-p_i$, and therefore generally gives suboptimal dependence on the sparsity parameters
$p_i$ (see Remark 2.4 in \cite{DHWZ_2025} for an example). The purpose of this paper is to derive concentration bounds that preserve this sparsity information.

\medskip

We now introduce the notation for the main theorem, adapted primarily from \cite{gotze2021concentration}. For a vector $x=(x_1,\ldots,x_n)^\top \in \mathbb R^n$, we write
$\|x\|_r=\left(\sum_i |x_i|^r\right)^{1/r}$ for the $\ell_{r}$ norm. 
For a random variable $\xi$, we write $\Vert \xi\Vert_{L_{r}}=(\mathbb{E}\vert \xi\vert^{r})^{1/r}$ for its $L_{r}$ norm.

Unless otherwise stated,  we denote by $C, C_{1}, c, c_{1},\ldots$ universal constants independent of any parameters and the dimension $n$. Similarly,  $C(\delta), c(\delta),\ldots $ denote constants depending only on the parameter $\delta$. Their values can change from line to line. We write $f\lesssim g$ if $f\le Cg$ for some universal constant $C$, and $f\lesssim_{\delta} g$ if $f\le C(\delta)g$ for some constant $C(\delta)$ depending only on $\delta$. We use $f\asymp g$ if $f\lesssim g$ and $g\lesssim f$, and similarly $f\asymp_{\delta} g$.

Let $[n]:=\{1,\ldots,n\}$. For $d\ge1 $, a multi-index
$ \mathbf i=(i_1,\ldots,i_d)\in [n]^d$ selects one coordinate in each of the $d$ directions.
A $d$-tensor is an array $A=(a_{\mathbf i})_{\mathbf i\in[n]^d}.$ 
For two $d$-tensors $A=(a_{\mathbf i})$ and $B=(b_{\mathbf i})$, their Hadamard product is
the $d$-tensor $A\circ B$ with entries
\[
(A\circ B)_{\mathbf i}=a_{\mathbf i}b_{\mathbf i}, \qquad \mathbf i\in[n]^d.
\]

Given an index set $I \subseteq [d]$, we write $\mathbf \mathbf{i}_I=(i_j)_{j\in I}.$
For a $d$-tensor $A$, the notation $A_{\mathbf \mathbf{i}_I}$ denotes the tensor obtained by fixing the coordinates in $I^c$ and allowing the coordinates in $ I$ to vary. Equivalently,
\(\mathbf \mathbf{i}_I\) records the free coordinates of the slice. For example, if $d=4$,
$I=\{2,4\}$, and $i_1=1,i_3=2$, then
$A_{\mathbf \mathbf{i}_I}=(a_{1\,i_2\,2\,i_4})_{i_2,i_4}.$ 
In particular, in the main theorem below, the expression
$ (\cdot)_{\mathbf i_{I^c}}$ denotes the slice obtained after fixing the coordinates indexed by $I$.

Let $\Pi_d$ be the set of all partitions of $[d]$.  More generally, for $I \subseteq [d]$, let
$ \Pi(I)$ be the set of all partitions of $I$. If
$ \mathcal J=\{J_1,\ldots,J_k\} \in \Pi(I) $, then $ |\mathcal J|= k$ denotes the number of blocks.
We use the convention that \(\Pi(\emptyset)=\{\emptyset\}\), and the corresponding tensor norm
of a scalar is its absolute value.

Given a partition $\mathcal J=\{J_1,\ldots,J_k\}\in\Pi_d$, we define
\[
\|A\|_{\mathcal J}:=
\sup\left\{
\sum_{\mathbf i\in[n]^d}
a_{\mathbf i}\prod_{\ell=1}^k x^{(\ell)}_{\mathbf i_{J_\ell}}: \sum_{\mathbf i_{J_\ell}}
\left(x^{(\ell)}_{\mathbf i_{J_\ell}}\right)^2\le 1,\;
\ell=1,\ldots,k
\right\}.
\]
The same definition applies to sliced tensors $A_{\mathbf \mathbf{i}_I}$ and partitions
$\mathcal J \in \Pi(I)$.

These norms increase as partitions become coarser. A partition $\mathcal I=\{I_1,\dots,I_\mu\}$ is finer than a partition $\mathcal J=\{J_1,\dots,J_\nu\}$ if  each block of $\mathcal J$ is a union of blocks from $\mathcal I$.  For example, $\{\{1,2\},\{3\},\{4\}\}$ is finer than $\{\{1,2,4\},\{3\}\}$. This gives the monotonicity
\[
\|A\|_{\mathrm{op}}=\|A\|_{\{\{1\},\dots,\{d\}\}}\le \|A\|_{\mathcal J}\le \|A\|_{\{[d]\}}=\|A\|_{\mathrm{HS}}=\Bigl(\sum_{\mathbf i}a_{\mathbf i}^2\Bigr)^{1/2}.
\]

For a multi-index $\mathbf i=(i_1,\ldots,i_d)\in [n]^d$, define
\[
\mathfrak p_{\mathbf i}^{(d)} := \sqrt{p_{j_1}\cdots p_{j_\nu}},
\]
where $j_1,\ldots,j_\nu$ are the distinct coordinate values appearing in
$ \mathbf i $. Thus repeated coordinates are counted only once. We write
\[
 \mathfrak{p}^{(d)} = \bigl(\mathfrak p_{\mathbf i}^{(d)}\bigr)_{\mathbf i\in[n]^d}
\]
for the corresponding $ d$-tensor. For $I \subseteq [d]$, define the partial sparsity weight
\[
\mathfrak p_I^{(d)}(\mathbf \mathbf{i}_I):= \sqrt{p_{k_1}\cdots p_{k_\mu}},
\]
where $ k_1,\ldots,k_\mu $ are the distinct coordinate values appearing in
$\mathbf \mathbf{i}_I=(i_j)_{j\in I}$. We use the convention $\mathfrak p_\emptyset^{(d)}=1.$ The role of the quotient
$ \mathfrak{p}^{(d)}/\mathfrak p_I^{(d)}$ will be illustrated after the theorem in the cubic case. 

We first state the general result. For $1\le d\le D$, $ I \subseteq [d]$, and
$\mathcal J\in\Pi([d]\setminus I)$, define
\[
M_{d,I,\mathcal J} := \max_{\mathbf \mathbf{i}_I} \left( \left\|\left(\mathbb E f^{(d)}(X)\circ  \mathfrak{p}^{(d)}
\right)_{\mathbf i_{I^c}}\right\|_{\mathcal J}\frac{1}{\mathfrak p_I^{(d)}(\mathbf \mathbf{i}_I)}\right).
\]
With the slicing convention above, the subscript $\mathbf i_{I^c}$ indicates that the coordinates in $I^c$ remain free, while the coordinates in $I$ are fixed and maximized over.

\begin{theorem}\label{Theo_main1}
Let $f:\mathbb R^n\to\mathbb R$ be a polynomial of degree $D$. 
Let $X_1,\cdots,X_n$ be a sequence of independent random variables such that $X_i=\delta_i\cdot\xi_i$, where the variables $\{\delta_i, \xi_i: i\in [n]\}$ are mutually independent, $\delta_{i}\sim \mathrm{Bernoulli}(p_i)$, and $\xi_i$ is an $\alpha$-sub-exponential random variable satisfying $\Vert \xi_i\Vert_{\Psi_\alpha}\le 1$.  Assume $0<\alpha\le 1$. Then, for any $r\ge 1$,
\[
\|f(X)-\mathbb E f(X)\|_{L_r} \lesssim_{D, \alpha}
\sum_{d=1}^D \sum_{I\subseteq[d]} \sum_{\mathcal J\in\Pi([d]\setminus I)}
r^{|\mathcal J|/2+|I|/\alpha} M_{d,I,\mathcal J}.
\]
\end{theorem}

\paragraph{The cubic case.}
We now illustrate Theorem~\ref{Theo_main1} in the first genuinely higher-order case. Since $D=2$ recovers the sparse Hanson-Wright bounds from \cite{DHWZ_2025}, the cubic case $D=3$ is the simplest setting where the new higher-order structure appears.

Assume, for this example, that $p_1=\cdots=p_n=p$, that the variables $X_i$ are centered,
and that
\[
f(x)=\frac1{3!}\sum_{i,j,k=1}^n B_{ijk}x_i x_j x_k,
\]
where $B$ is symmetric and diagonal-free, i.e., $B_{ijk}=0$ whenever $|\{i,j,k\}|<3$. Then $f^{(3)}\equiv B$, while $ \mathbb E f^{(1)}(X)=\mathbb E f^{(2)}(X)=0.$  For any index vector $\mathbf{i}$ corresponding to a nonzero entry of $B$, we have 
$\frac{\mathfrak p^{(3)}_{\mathbf i}}{\mathfrak p_I^{(3)}(\mathbf \mathbf{i}_I)}=p^{(3-|I|)/2}.$

Theorem~\ref{Theo_main1} yields that, for every $r\ge1$, 
\begin{align*}
\|f(X)-\mathbb E f(X)\|_{L_r} \lesssim_{\alpha}\;& p^{3/2} \left[ r^{1/2} \|B\|_{\{[3]\}}
+ r \sum_{\substack{\mathcal J \in \Pi([3]) \\ |\mathcal J|=2 }} \|B\|_{\mathcal J}
+ r^{3/2} \|B\|_{\{\{1\},\{2\},\{3\}\}} \right] \\
&+ p \sum_{\substack{I\subseteq[3]\\ |I|=1}} \left[ r^{1/2+1/\alpha} \max_{\mathbf \mathbf{i}_I} \|B_{\mathbf i_{I^c}}\|_{\{I^c\}}
+ r^{1+1/\alpha} \max_{\mathbf \mathbf{i}_I}
\|B_{\mathbf i_{I^c}}\|_{\{\{j\}:j\in I^c\}} \right] \\
&+ \sqrt p\, r^{1/2+2/\alpha}
\sum_{\substack{I\subseteq[3]\\ |I|=2}} \max_{\mathbf \mathbf{i}_I} \|B_{\mathbf i_{I^c}}\|_{\{I^c\}} + r^{3/\alpha} \max_{\mathbf i_{[3]}} |B_{i_1 i_2 i_3}|.
\end{align*}
To interpret these partition norms, note that \[\|B\|_{\{[3]\}} = \|B\|_{\rm HS},\qquad \|B\|_{\{\{1\},\{2\},\{3\}\}} = \|B\|_{\rm op},\] and the middle terms $\|B\|_{\mathcal J}$ for $|\mathcal J|=2$ interpolate between them. When coordinates are fixed to form slices $\|B_{\mathbf i_{I^c}}\|$, different structures emerge. If $|I|=1$, fixing one coordinate yields matrix slices with norms $\|B_{\mathbf i_{I^c}}\|_{\{I^c\}}$ (Hilbert-Schmidt norm) and $\|B_{\mathbf i_{I^c}}\|_{\{\{j\}:j\in I^c\}}$ (operator norm). If $|I|=2$, fixing two coordinates yields vector fibers with $\ell_2$ norm $\|B_{\mathbf i_{I^c}}\|_{\{I^c\}}$. If $|I|=3$, fixing all coordinates yields the maximum entry $\|B\|_{\max}=\max_{i_1, i_2, i_3} |B_{i_1 i_2 i_3}|$.

This display illustrates the structure of the general theorem. If no coordinate is fixed, all three Bernoulli selectors contribute to the averaging scale, yielding the factor $p^{3/2}$. If one coordinate is fixed, the remaining matrix slice carries the factor $p$. If two coordinates are fixed, the remaining vector fiber carries the factor $\sqrt p$. If all three coordinates are fixed,
there is no uniform sparsity gain in the high-moment regime, and the corresponding term is
controlled by $\|B\|_{\max}$. Each fixed coordinate costs a heavy-tail factor $r^{1/\alpha}$, while the remaining free coordinates are measured by partition norms.

\medskip

Theorem \ref{Theo_main1} and the standard moment-tail relationship (see, e.g., Proposition 3.3 in \cite{gotze2021concentration}) immediately yield the following exponential tail bound for $f(X)$. We omit the proof via Markov's inequality as it is standard.

\begin{coro}\label{Coro_main_tail}
Under the assumptions of Theorem~\ref{Theo_main1}, there exist constants
$C_{D,\alpha},c_{D,\alpha}>0$ such that, for any $ t >0$,
\[
\mathbb P\left\{
|f(X)-\mathbb E f(X)|\ge t
\right\} \le C_{D,\alpha} \exp\left[ -c_{D,\alpha}
\min_{\substack{1\le d\le D,\; I\subseteq[d],\\
\mathcal J\in\Pi([d]\setminus I),\; M_{d,I,\mathcal J}>0}}
\left( \frac{t}{M_{d,I,\mathcal J}} \right)^{ \frac{1}{|\mathcal J|/2+|I|/\alpha} } \right].
\]
If all $M_{d,I,\mathcal J}$ vanish, then $f(X)-\mathbb E f(X)=0$ almost surely and the
bound is trivial.
\end{coro}

\subsection{Discussion of Theorem~\ref{Theo_main1}}

When $p_1=\cdots=p_n=1$, all sparsity weights are equal to one. Hence
Theorem~\ref{Theo_main1} reduces to the polynomial concentration inequality of
G{\" o}tze, Sambale, and Sinulis \cite{gotze2021concentration} for $\alpha$-sub-exponential inputs. In this sense,
Theorem~\ref{Theo_main1} is a sparse analogue of their result.

The case $D=2$ recovers the sparse Hanson-Wright inequality proved in
\cite{DHWZ_2025}. More precisely, for $\mathbb{E}\xi_i = 0$ and symmetric coefficients $a_{ij}=a_{ji}$, Theorem \ref{Theo_main1} yields a four-term estimate for $f(X)=\sum_{i,j} a_{ij}X_i X_j$:
\begin{align}\label{eq:quadratic_degeneration}
\left\| f(X)-\mathbb E f(X)\right\|_{L_r}
&\lesssim_{\alpha} r^{2/\alpha}\max_{i,j}|a_{ij}| + r^{1/2+1/\alpha} \max_i \Big(\sum_{j=1}^n a_{ij}^2p_j\Big)^{1/2} \nonumber \\
&\quad + r\left\|(a_{ij}\sqrt{p_ip_j})_{i,j}\right\|_{\mathrm{op}} + \sqrt{r} \Big(\sum_{i=1}^n a_{ii}^2p_i + \sum_{i\ne j}a_{ij}^2p_ip_j\Big)^{1/2}.
\end{align}
This illustrates the sparse Hanson-Wright phenomenon: averaged coordinates contribute powers of $p_i^{1/2}$, while coordinates controlled through high moments do not carry a uniform sparsity gain (see \cite{He_arXiv_HW,DHWZ_2025} for detailed discussions). 

For general $D$, the same mechanism is encoded by the quotient
\[
{ \mathfrak{p}^{(d)}}/{\mathfrak p_I^{(d)}}.
\]
The numerator $ \mathfrak{p}^{(d)}$ records the Bernoulli weights for all distinct
coordinates in a derivative tensor. The denominator $\mathfrak p_I^{(d)}$ removes the sparsity
gain from the coordinates fixed in the high-moment part. This removal is necessary: for instance,
a single monomial like $X_1^D=\delta_1\xi_1^D$ has $ L_r$-norm of order $ r^{D/\alpha}$
in the regime $ r\gtrsim D\log(1/p_1) $, so no positive power of $ p_1$ can multiply the
corresponding extreme-coordinate term uniformly in $r$.

The appearance of all derivative orders $1\le d\le D$ is also unavoidable for non-centered
polynomials. For example, if $X_i=\delta_i\sim\mathrm{Bernoulli}(p)$, $S=\sum_{i=1}^n X_i$, and
$ f(x)=(\sum_i x_i)^D$, then the leading fluctuation of $S^D$ is of order $(np)^{D-1}\sqrt{np}.$ 
This scale is captured by the first derivative $\mathbb E f^{(1)}(X) $, not by the top derivative
alone. Thus the lower-order derivative tensors in Theorem~\ref{Theo_main1} are a structural
feature of polynomial concentration for non-centered inputs, already present in the dense results
of Adamczak-Wolff \cite{AW15} and G\" otze-Sambale-Sinulis \cite{gotze2021concentration}.

Finally, we comment on the proof strategy. Our proof follows the polynomial concentration framework from the dense setting (see \cite{gotze2021concentration} and reference therein): one expands the polynomial into centered monomial chaoses, applies decoupling and symmetrization, reduces the variables to Weibull-type
models through a contraction principle, and then invokes moment estimates for chaoses with
log-convex tails. The new difficulty in the sparse setting is that the Bernoulli selectors must be
kept throughout this reduction. Their contribution depends on how many distinct coordinates
remain averaged and how many coordinates are fixed in the high-moment regime. This is exactly
what produces the distinct-coordinate weight $ \mathfrak{p}^{(d)}$ and its quotient by
$\mathfrak p_I^{(d)}$. The key technical contribution is to propagate these weights through the
partition-norm estimates and then convert the resulting coefficient tensors back into the expected
derivative tensors $\mathbb E f^{(d)}(X)$.

\section{Applications: Subspace Distance and Smallest Singular Values}\label{sec:applications}

In this section we illustrate Theorem~\ref{Theo_main1} through a geometric concentration result for sparse simple random tensors. Let
\[
     X = X^{(1)}\otimes\cdots\otimes X^{(d)}\in \mathbb R^{n^d},
\]
where the vectors \(X^{(1)},\ldots,X^{(d)}\in\mathbb R^n\) have independent sparse \(\alpha\)-sub-exponential entries. 

A natural problem is to understand the distance from \( X\) to a fixed subspace \(H\subseteq \mathbb R^{n^d}\). Writing \(P_{H^\perp}\) for the orthogonal projection onto \(H^\perp\), we have
\[
    \operatorname{dist}( X,H)=\|P_{H^\perp} X\|_2.
\]
For $d=1$, this is a Lipschitz function of the underlying independent entries. For $d\ge 2$,
however, the map
\[
    (X^{(1)},\ldots,X^{(d)}) \mapsto \|P_{H^\perp}(X^{(1)}\otimes\cdots\otimes X^{(d)})\|_2
\]
is no longer globally Lipschitz in the scalar inputs. Its square is a polynomial of degree $2d$. More generally, for a deterministic matrix $ B\in \mathbb R^{n^d\times n^d}$,
\[
    \|BX\|_2^2 =  X^\top B^\top B X .
\]
Thus Theorem~\ref{Theo_main1} applies naturally to this problem and yields a deviation inequality which keeps track of the sparsity parameter.

%Distance estimates of this form are useful in random matrix problems. In particular, by the standard leave-one-out argument and the negative second moment identity, they imply lower bounds for the smallest singular value of matrices whose columns are independent sparse simple random tensors. We present this consequence as a corollary after the distance theorem.
The following theorem gives the resulting concentration estimate. Throughout this section, the
tensor order $d$ is fixed, and constants may depend on $d$ and $\alpha$.

\begin{theorem}[Distance theorem for sparse simple random tensors]
\label{thm:distance}
Let \(0<\alpha\le 1\), \(p\in(0,1]\), and let
\[
     X=X^{(1)}\otimes\cdots\otimes X^{(d)}\in\mathbb R^{n^d},
\]
where \(X_i^{(k)}=\delta_i^{(k)}\xi_i^{(k)}\), the variables \(\delta_i^{(k)}\) are i.i.d. Bernoulli\((p)\), and the variables \(\xi_i^{(k)}\) are independent, mean zero, variance one, and satisfy  $\|\xi_i^{(k)}\|_{\Psi_\alpha}\le 1$. All variables are assumed to be mutually independent.
Set
\[
    \Lambda_{n,p,d}:=\max_{0\le s\le d-1}(np)^{s/2}.
\]
Then for every deterministic nonzero matrix \(B\in\mathbb R^{n^d\times n^d}\) and every \(t\ge 0\),
\begin{align}
&\mathbb P\left\{
\left|\|B X\|_2-p^{d/2}\|B\|_{\mathrm{HS}}\right|\ge t
\right\}\nonumber\\
&\qquad\le 2\exp\left[
-c_{\alpha,d}
\min\left\{
\frac{t^2}{\Lambda_{n,p,d}^2\|B\|_{\mathrm{op}}^2},
\left(\frac{t}{n^{(d-1)/2}\|B\|_{\mathrm{op}}}\right)^\alpha,
\left(\frac{t}{\|B\|_{\mathrm{op}}}\right)^{\alpha/d}
\right\}
\right].
\end{align}
In particular, if \(np\ge 1\), then
\[
    \Lambda_{n,p,d}=(np)^{(d-1)/2}.
\]
\end{theorem}

Taking $B=P_{H^\perp}$ in Theorem~\ref{thm:distance} gives the following distance estimate.
If $H\subseteq\mathbb R^{n^d}$ is a fixed subspace and $k:=\dim (H^\perp),$ 
then $\|P_{H^\perp}\|_{\mathrm{op}}=1$ and  $\|P_{H^\perp}\|_{\mathrm{HS}}=\sqrt{k}.$ 
Consequently, for every $t \ge 0$,
\begin{align*}
&\mathbb P\left\{
\left|\operatorname{dist}(X,H)-p^{d/2}\sqrt{k}
\right|\ge t\right\} \\
&\qquad\le 2\exp\left[ -c_{\alpha,d} \min\left\{\frac{t^2}{\Lambda_{n,p,d}^2},
\left(\frac{t}{n^{(d-1)/2}}\right)^\alpha, t^{\alpha/d}\right\} \right].
\end{align*}

The distance estimate of the preceding form serves as a core tool for studying random matrices with tensor structures. Indeed, by the standard leave-one-out argument and the negative second moment
identity, lower bounds on
\[
 \operatorname{dist}(\mathbf X_j,H_j), \qquad H_j:=\operatorname{span}\{\mathbf X_i:i\ne j\},
\]
imply lower bounds for the smallest singular value of the column matrix.

For $d=1$, this problem reduces to the smallest singular value of an ordinary rectangular random matrix, which has been extensively investigated. For independent sub-Gaussian entries, Rudelson and Vershynin \cite{RV2009CPAM} proved lower bounds of the optimal order. In the sparse regime, G{\"o}tze and Tikhomirov \cite{GT2023EJP} demonstrated $\sigma_{\min}(X) \ge \sqrt{mp}$ under certain moment conditions. For higher-order tensors ($d \ge 2$), Vershynin \cite{vershynin2020concentration} established that $\sigma_{\min}(X) \ge \sqrt{\epsilon}$ for dense sub-Gaussian tensor matrices when $m = (1-\epsilon)n^d$. 

As a direct consequence of Theorem 2, we extend these frameworks simultaneously to the sparse and heavy-tailed setting, providing a non-trivial lower bound for sparse simple random tensor matrices.

\begin{coro}[Smallest singular value of sparse simple tensor matrices]
\label{cor:smin}
Let \(\mathbf X_1,\ldots,\mathbf X_m\) be independent copies of
\[
     X=X^{(1)}\otimes\cdots\otimes X^{(d)}\in\mathbb R^{n^d}
\]
as in Theorem~\ref{thm:distance}, and let
\[
    \mathbf X=(\mathbf X_1,\ldots,\mathbf X_m)\in\mathbb R^{n^d\times m}.
\]
Suppose that
\[
    \frac{c_1(d\log n)^{2/\alpha}}{np^d}\le \varepsilon<1
\]
and
\[
    m\le (1-\varepsilon)n^d.
\]
Then
\[
\mathbb P\left\{
\sigma_{\min}(\mathbf X)>
\frac12 p^{d/2}\sqrt{\frac{\varepsilon}{1-\varepsilon}}
\right\}
\ge
1-2\exp\left[-c_2(p^d\varepsilon n)^{\alpha/2}\right],
\]
where \(c_1,c_2>0\) depend only on \(\alpha\) and \(d\).
\end{coro}

\paragraph{Remarks on the applications.} Theorem~\ref{thm:distance} applies sparse polynomial concentration to quadratic forms in simple random tensors. If $X = X^{(1)}\otimes\cdots\otimes X^{(d)} \in \mathbb R^{n^d}$, then $\|BX\|_2^2 = X^\top B^\top BX$ is a degree-$2d$ polynomial in the underlying independent entries. The proof in Section~\ref{subsec:tensor_concentration} establishes a more general moment estimate (see Theorem~\ref{thm:tensor_poly_concentration}) for such tensor quadratic forms, but we defer that technical statement to the proof section since it requires additional partial-trace notation. For geometric applications, Theorem~\ref{thm:distance} provides the main tool.

Theorem~\ref{thm:distance} gives a sparse, heavy-tailed analogue of tensor Hanson--Wright concentration. The closest dense result is the tensor Hanson-Wright inequality of Bamberger, Krahmer and Ward~\cite{BKW22} for quadratic forms in simple random tensors with independent centered sub-Gaussian entries. Our setting differs: the entries are sparse and only $\alpha$-sub-exponential with $0<\alpha\le1$. Accordingly, the typical distance scale becomes $p^{d/2}\|B\|_{\rm HS}$, and the leading fluctuation contains the sparse factor $\Lambda_{n,p,d}=\max_{0\le s\le d-1}(np)^{s/2}$. Thus our result complements rather than improves the dense sub-Gaussian theory.

The tail bound in Theorem~\ref{thm:distance} has three regimes corresponding to \[ t^2/(\Lambda_{n,p,d}^2\|B\|_{\rm op}^2), \quad (t/(n^{(d-1)/2}\|B\|_{\rm op}))^\alpha, \quad (t/\|B\|_{\rm op})^{\alpha/d}.\] The first captures Gaussian-type fluctuations, while the second and third reflect heavy-tailed behavior from products of $\alpha$-sub-exponential variables. These additional regimes are natural for degree-$d$ tensor products but absent in purely sub-Gaussian linear models.

The constants in our results depend on the tensor order $d$, so they are most useful when $d$ is fixed. This differs from dense sub-Gaussian tensor results like Vershynin~\cite{vershynin2020concentration}, where the dependence on $d$ is tracked explicitly. The advantage of our approach is its applicability to sparse inputs and the heavy-tailed regime $0<\alpha\le1$.

Corollary~\ref{cor:smin} should be viewed as a consequence of general polynomial concentration rather than an optimal smallest-singular-value theorem. For ordinary matrices ($d=1$), sharper bounds exist through specialized methods. The value of Corollary~\ref{cor:smin} is that the same distance argument works for tensor-product columns with sparse $\alpha$-sub-exponential entries. The probability estimate becomes meaningful once $p^d\varepsilon n$ is at least logarithmic in $n$, which motivates the assumption $\varepsilon\ge c_1(d\log n)^{2/\alpha}/(np^d)$.

\section{Proofs of our main results}\label{sec:proofs}
\subsection{Proof of Theorem \ref{Theo_main1}}
To begin, we introduce the notation and auxiliary lemmas that will be used throughout this section. For any multiindex $\mathbf{i}$, denote $|\mathbf{i}|=\sum_{j}i_j$.
Write
\[
[n]^{\underline{d}}:=\bigl\{\mathbf{i}\in [n]^{d}:i_{1},\dots ,i_{d}\text{ are pairwise distinct}\bigr\}
\] and 
\begin{align}
    I_{m, d}:=\{(i_1, \cdots, i_m)\in \mathbb{N}_{+}^m: |\mathbf{i}|=d \}.
\end{align}
For two \(d\)-tensors \(A=(a_{\mathbf{i}})\) and \(B=(b_{\mathbf{i}})\), their (Frobenius) inner product is
\[
\langle A,B\rangle:=\sum_{\mathbf{i}\in [n]^{d}}a_{\mathbf{i}}\,b_{\mathbf{i}}.
\]

The first property we shall introduce is the following contraction principle.
\begin{lemma}[Lemma 3.1 in \cite{gotze2021concentration}]\label{Lem_contration}
   For every \(k\in\mathbb{N}\), \(\alpha>0\) and \(r\ge 1\) the following holds.  
Let \(Y_{1},\dots ,Y_{n}\) be independent, symmetric random variables satisfying
$
\|Y_{i}\|_{\Psi_{\alpha/k}}\le M.
$
Then
\[
\Bigl\|\sum_{i=1}^{n}a_{i}Y_{i}\Bigr\|_{L_{r}}
\ \lesssim_{\alpha,k,M}\ 
\Bigl\|\sum_{i=1}^{n}a_{i}\,w_{i,1}\cdots w_{i,k}\Bigr\|_{L_{r}},
\]
where the \(w_{i,j}\) are i.i.d. symmetric Weibull random variables with shape parameter \(\alpha\).
\end{lemma}
 Recall that a random variable \(\xi\) is said to have a log-convex tail if the function  
\[
t \mapsto \log \mathbb{P}\{|\xi| \ge t\}
\]  
is convex for \(t \ge 0\).  
We will also need the following properties which bounds  the \(L_r\)-norms of multilinear forms in random variables with log-convex tails.

\begin{lemma}[Theorem 1 in \cite{KL15}]\label{Lem_Lrbound}
    Let $(X_{i}^{(j)})_{i\le n, j\le d}$ be independent symmetric r.v.'s with log-convex tails such that $\mathbb{E}\vert X_{i}^{(j)}\vert^2=1$ for all $i , j$. For any multiindexed matrix $(a_\mathbf{i})_{\mathbf{i}\in [n]^d}$ and any $r\ge 2$ we have
    \begin{align}
        \Big\Vert \sum_{\mathbf{i}}a_{\mathbf{i}}X_{i_1}^{(1)}\cdots X_{i_d}^{(d)} \Big\Vert_{L_r}\asymp_{d} \sum_{I\subseteq [d]}\sum_{\mathcal{J}\in \Pi(I^c)}r^{|\mathcal{J}|/2} \big(\sum_{\mathbf{i}_{I}}\Vert (a_{\mathbf{i}})_{\mathbf{i}_{I^c}}\Vert_{\mathcal{J}}^r\prod_{j\in I}\Vert X_{i_j}^{(j)}\Vert_{L_r}^{r} \big)^{1/r}.
    \end{align}
\end{lemma}

To state the next property we fix some notation.  
For any subset \(S\subseteq[n]^{d}\) let \(\mathbf{1}_{S}\) be the indicator tensor whose entry at position \(\mathbf{i}\in[n]^{d}\) is \(1\) if \(\mathbf{i}\in S\) and \(0\) otherwise.  
Given a partition \(\mathcal{K}=\{K_{1},\dots,K_{\nu}\}\) of \([d]\), we set

\[
L(\mathcal{K})=\bigl\{\mathbf{i}\in[n]^{d}:i_{k}=i_{l}\;\Longleftrightarrow\;\exists\,j\text{ such that }k,l\in K_{j}\bigr\}.
\]
In words, \(L(\mathcal{K})\) collects all multi-indices whose level-set partition is exactly \(\mathcal{K}\).

\begin{lemma}[Lemma 3.4 in \cite{gotze2021concentration}]\label{Lem_matrixnorm_comparasion}
    Let $A=(a_\mathbf{i})_{\mathbf{i}\in[n]^d}$ be a $d$-tensor, $I\subseteq [d]$ and $\mathbf{i}_I\in [n]^{|I|}$ fixed. Then for any $\mathcal{K}\in \Pi_d$ and $\mathcal{J}\in \Pi(I^c)$,
    \begin{align}
        \Vert (A\circ \mathbf{1}_{L(\mathcal{K})})_{\mathbf{i}_{I^c}}\Vert_{\mathcal{J}}\le 2^{|\mathcal{K}|(|\mathcal{K}|-1)/2}\Vert A_{\mathbf{i}_{I^c}}\Vert_{\mathcal{J}}.
    \end{align}
\end{lemma}

\begin{proof}[Proof of Theorem \ref{Theo_main1}]
Without loss of generality, we can represent any order $D$ polynomial $f$ as
\begin{align}
    f(x)=\sum_{d=1}^{D}\sum_{\nu=1}^{d}\sum_{\mathbf{k}\in I_{\nu, d}}\sum_{\mathbf{i}\in [n]^{\underline{\nu}}}c^{(d)}_{(i_1, k_1), \cdots, (i_{\nu}, k_{\nu})}x_{i_1}^{k_1}x_{i_2}^{k_2}\cdots x_{i_\nu}^{k_\nu}+c_0,
\end{align}
where the constant satisfies that $c^{(d)}_{(i_1, k_1), \cdots, (i_{\nu}, k_{\nu})}=c^{(d)}_{(i_{\pi_1}, k_{\pi_1}), \cdots, (i_{\pi_\nu}, k_{\pi_\nu})}$ for any permutation $(\pi_{1}, \cdots, \pi_{\nu})$ of $[\nu]$. Note that
\begin{align}
    &f(X)-\mathbb{E}f(X)\nonumber\\
    =&\sum_{d=1}^D\sum_{\nu=1}^d\sum_{\mathbf{k}\in I_{\nu, d}}\sum_{\mathbf{i}\in[n]^{\underline{\nu}}}c^{(d)}_{(i_1, k_1), \cdots, (i_{\nu}, k_{\nu})} \prod_{j\in [\nu]}(X_{i_j}^{k_j}-\mathbb{E}X_{i_j}^{k_j}+\mathbb{E}X_{i_j}^{k_j})\nonumber\\
    =&\sum_{d=1}^D\sum_{\nu=1}^d\sum_{\mathbf{k}\in I_{\nu, d}}\sum_{\mathbf{i}\in[n]^{\underline{\nu}}}c^{(d)}_{(i_1, k_1), \cdots, (i_{\nu}, k_{\nu})}\left( \sum_{\emptyset\neq J\subseteq [\nu]}\prod_{j\in J}(X_{i_j}^{k_j}-\mathbb{E}X_{i_j}^{k_j})\prod_{j\notin J}\mathbb{E}X_{i_j}^{k_j}\right).
\end{align}
Hence, by rearranging the terms and using the triangle inequality, we have 
\begin{align}
    \big\vert f(X)-\mathbb{E} f(X) \big\vert\le \sum_{d=1}^D\sum_{\nu=1}^d\sum_{\mathbf{k}\in I_{\nu, d}}\Big\vert\sum_{\mathbf{i}\in [n]^{\underline{\nu}}}a_{\mathbf{i}}^{\mathbf{k}}(X_{i_1}^{k_1}-\mathbb{E}X_{i_1}^{k_1})\cdots (X_{i_\nu}^{k_\nu}-\mathbb{E}X_{i_\nu}^{k_\nu})\Big\vert,
\end{align}
    where 
    \begin{align}
        a_{\mathbf{i}}^{\mathbf{k}}=a_{i_1,\cdots, i_\nu}^{(k_1, \cdots, k_\nu)}:=\sum_{m=\nu}^{D}\sum_{\substack{
    k_{\nu +1}, \cdots, k_m>0\\[2pt]    
    k_1+\cdots+k_m\le D}}\sum_{\substack{
    i_{\nu +1}, \cdots, i_m\\[2pt]    
    (i_1,\cdots,i_m)\in [n]^{\underline{m}}}}\binom{m}{\nu}c^{(k_1+\cdots+k_{m})}_{(i_1, k_1),\cdots, (i_m, k_m)}\prod_{\beta=\nu+1}^m\mathbb{E}X_{i_{\beta}}^{k_{\beta}}.
    \end{align}
    Indeed, in the rearrangement of the sum, for a given degree $m$ with $\nu$ many terms of $(X_i^k-\E X_i^k)$, we will just consider the pattern $\prod_{j=1}^\nu(X_{i_j}^{k_j}-\mathbb{E}X_{i_j}^{k_j})\prod_{j=\nu+1}^m\mathbb{E}X_{i_j}^{k_j}$. Then for given $m,\nu$, there will be $\binom m\nu$ many products corresponds to this pattern in the sum due to symmetry.

    Let $X^{(1)}, \cdots, X^{(D)}$ be independent copies of the random vector $X$. Let $(\varepsilon_{i}^{(j)}), 1\le i\le n, 1\le j\le D$ be a sequence of i.i.d. Rademacher variables independent of $X^{(1)}, \cdots X^{(D)}$. By the standard decoupling inequalities (see Theorem 3.1.1 in \cite{kwapien1987decoupling}) and symmetrization technique (see Section 6.1 in \cite{Ledoux_probability_Banach}), we have for $r\ge 1$
    \begin{align}
        &\big\Vert f(X)-\mathbb{E} f(X) \big\Vert_{L_r}\nonumber\\
        \le& \sum_{d=1}^D\sum_{\nu=1}^d\sum_{\mathbf{k}\in I_{\nu, d}}\Big\Vert\sum_{\mathbf{i}\in [n]^{\underline{\nu}}}a_{\mathbf{i}}^{\mathbf{k}}((X_{i_1}^{(1)})^{k_1}-\mathbb{E}(X_{i_1}^{(1)})^{k_1})\cdots ((X_{i_\nu}^{(\nu)})^{k_\nu}-\mathbb{E}(X_{i_\nu}^{(\nu)})^{k_\nu})\Big\Vert_{L_r}\nonumber\\
        \lesssim& \sum_{d=1}^D\sum_{\nu=1}^d\sum_{\mathbf{k}\in I_{\nu, d}}\Big\Vert\sum_{\mathbf{i}\in [n]^{\underline{\nu}}}a_{\mathbf{i}}^{\mathbf{k}}\varepsilon_{i_1}^{(1)}(X_{i_1}^{(1)})^{k_1}\cdots \varepsilon_{i_\nu}^{(\nu)}(X_{i_\nu}^{(\nu)})^{k_\nu}\Big\Vert_{L_r}\nonumber\\
        \lesssim& \sum_{d=1}^D\sum_{\nu=1}^d\sum_{\mathbf{k}\in I_{\nu, d}}\Big\Vert\sum_{\mathbf{i}\in [n]^{\underline{\nu}}}a_{\mathbf{i}}^{\mathbf{k}}\delta_{i_1}^{(1)}\varepsilon_{i_1}^{(1)}(\xi_{i_1}^{(1)})^{k_1}\cdots \delta_{i_\nu}^{(\nu)}\varepsilon_{i_\nu}^{(\nu)}(\xi_{i_\nu}^{(\nu)})^{k_\nu}\Big\Vert_{L_r}.
    \end{align}
    Due to that $\Vert (\xi_{i}^{(j)})^k\Vert_{\Psi_{\alpha/k}}=\Vert\xi_{i}^{(j)}\Vert_{\Psi_{\alpha}}^k\le 1, j=1,\cdots, D $, we have by an iteration of  Lemma \ref{Lem_contration}
    \begin{align}
        &\big\Vert f(X)-\mathbb{E} f(X) \big\Vert_{L_r}\nonumber\\
        \lesssim &\sum_{d=1}^D\sum_{\nu=1}^d\sum_{\mathbf{k}\in I_{\nu, d}}\Big\Vert\sum_{\mathbf{i}\in [n]^{\underline{\nu}}}a_{\mathbf{i}}^{\mathbf{k}}(\delta_{i_1}^{(1)}w_{i_1, 1}^{(1)}\cdots w_{i_1, k_1}^{(1)})\cdots (\delta_{i_\nu}^{(\nu)}w_{i_\nu, 1}^{(\nu)}, \cdots w_{i_\nu, k_\nu}^{(\nu)})\Big\Vert_{L_r}, \label{Eq_proof_theo1_1}
    \end{align}
  where $\{w_{i, k}^{(j)}, 1\le i, j, k\le D\}$ is a sequence of i.i.d. symmetric Weibull variables with shape parameter $\alpha$, and the extra factor of $\eps$ will not change its distribution.  

  We next define a $d$-tensor $A_d=(a^{(d)}_{\mathbf{i}})_{\mathbf{i}\in [n]^d}$ as follows. Given any $\mathbf{i}=(i_1, \cdots, i_d)$, there are distinct elements $j_1, \cdots, j_\nu$ with each $j_l$ appearing exactly $k_l$ times in $\mathbf{i}$. Then, we set $a^{(d)}_{i_1,\cdots, i_d}:=a_{j_1,\cdots, j_\nu}^{(k_1,\cdots, k_\nu)}$ (as $a_{\mathbf j}^{\mathbf k}$ is invariant under permutation on $[\nu]$). Obviously, we have $k_1+\cdots+k_\nu=d$. For any $\mathbf{k}=(k_1, \cdots, k_\nu)\in I_{\nu, d}$, denote by $\mathcal{K}(\mathbf{k})=\{K_1, \cdots, K_\nu\}\in \Pi_d$ such that $K_l=\{\sum_{i=1}^{l-1}k_i+1, \sum_{i=1}^{l-1}k_i+2, \cdots, \sum_{i=1}^l k_i\}, l=1, \cdots, \nu$. For convenience, given a $\mathbf{k}\in I_{\nu, d}$, let
  \begin{align}
      &Y^{(1)}=\big(w_{i, 1}^{(1)}\cdot\delta_{i}^{(1)}\big)_{i\le n}, Y^{(2)}=\big(w_{i, 2}^{(1)}\big)_{i\le n},\cdots, Y^{(k_1)}=\big(w_{i, k_1}^{(1)}\big)_{i\le n}, \nonumber\\
      &Y^{(k_1+1)}=\big(w_{i, 1}^{(2)}\cdot\delta_{i}^{(2)}\big)_{i\le n},Y^{(k_1+2)}=\big(w_{i, 2}^{(2)}\big)_{i\le n},\cdots, Y^{(k_1+k_2)}=\big(w_{i, k_2}^{(2)}\big)_{i\le n},\\
      &\cdots\\
      &Y^{(\sum_{j=1}^{\nu-1}k_j+1)}=\big(w_{i, 1}^{(\nu)}\cdot\delta_{i}^{(\nu)}\big)_{i\le n},Y^{(\sum_{j=1}^{\nu-1}k_j+2)}=\big(w_{i, 2}^{(\nu)}\big)_{i\le n}\cdots, Y^{(d)}=\big(w_{i, k_\nu}^{(\nu)} \big)_{i\le n}.
  \end{align}
Note that
$$		-\log \mathbb{P}\{\vert\delta^{(1)}_{i_1}w^{(1)}_{i_1, 1}\vert\ge t \} = 
		\begin{cases}
			0, & t = 0 \\
			t^{\alpha}-\log p_{i}, & t > 0
		\end{cases}
$$
is a concave function for $t\ge 0$.  Since Lemma~\ref{Lem_Lrbound}
is stated under the normalization $\mathbb E|X_i^{(j)}|^2=1$, we apply it
below to the normalized variables $Y_i^{(j)}/\|Y_i^{(j)}\|_2$; the deterministic
factors $\|Y_i^{(j)}\|_2$ are then absorbed into the coefficient tensor. Hence, we can further bound \eqref{Eq_proof_theo1_1} by Lemma \ref{Lem_Lrbound}
\begin{align}
    &\sum_{d=1}^D\sum_{\nu=1}^d\sum_{\mathbf{k}\in I_{\nu, d}}\Big\Vert\Big\langle A_d\circ \mathbf{1}_{L(\mathcal{K}(\mathbf{k}))}, \otimes_{j=1}^{d}Y^{(j)}  \Big\rangle\Big\Vert_{L_r}\nonumber\\
    \lesssim_d&\sum_{d=1}^D\sum_{\nu=1}^d\sum_{\mathbf{k}\in I_{\nu, d}}\sum_{I\subseteq [d]}\sum_{\mathcal{J}\in \Pi([d]\backslash I)}r^{|\mathcal{J}|/2}\Big( \sum_{\mathbf{i}_{I}}\Big(\big\Vert \big(A_d\circ \mathfrak{p}^{(d)}\circ \mathbf{1}_{L(\mathcal{K}(\mathbf{k}))}\big)_{\mathbf{i}_{I^c}}\big\Vert^r_{\mathcal{J}}\nonumber\\
    &\times\prod_{j\in I}\Vert Y_{i_j}^{(j)}\Vert_{L_r}^r\prod_{j\in I\cap S(\mathbf{k}) } \frac{1}{\sqrt{p^r_{i_j}}}\Big) \Big)^{1/r},
\end{align}
where $\mathfrak{p}^{(d)}$ is defined as in Theorem \ref{Theo_main1}, and $S(\mathbf{k})$ is defined as follows: for a vector \(\mathbf k=(k_1,\dots,k_\nu)\) we set \(S(\mathbf k)=\{1,k_1+1,k_1+k_2+1,\dots,k_1+\cdots+k_{\nu-1}+1\}\). 
\allowdisplaybreaks

Let $\tilde{Y}^{(1)}, \cdots, \tilde{Y}^{(d)}\overset{\text{i.i.d.}}{\sim}Y^{(2)}$. Then, by Young's inequality, we have for $r\ge 3$
\begin{align}
    &\Big( \sum_{\mathbf{i}_{I}}\big\Vert \big(A_d\circ \mathfrak{p}^{(d)}\circ \mathbf{1}_{L(\mathcal{K}(\mathbf{k}))}\big)_{\mathbf{i}_{I^c}}\big\Vert^r_{\mathcal{J}}\prod_{j\in I}\Vert Y_{i_j}^{(j)}\Vert_{L_r}^r\prod_{j\in I\cap S(\mathbf{k}) } \frac{1}{\sqrt{p^r_{i_j}}} \Big)^{1/r}\nonumber\\
    =&\Big( \sum_{\mathbf{i}_{I}}\big\Vert \big(A_d\circ \mathfrak{p}^{(d)}\circ \mathbf{1}_{L(\mathcal{K}(\mathbf{k}))}\big)_{\mathbf{i}_{I^c}}\big\Vert^r_{\mathcal{J}}\prod_{j\in I}\Vert \tilde{Y}_{i_j}^{(j)}\Vert_{L_r}^r\prod_{j\in I\cap S(\mathbf{k}) } \frac{p_{i_j }}{\sqrt{p^r_{i_j}}} \Big)^{1/r}\nonumber\\
    \lesssim_d& \Big(e^{2r/(r-2)}\max_{\mathbf{i}_I}\big\Vert\big(A_d\circ \mathfrak{p}^{(d)}\circ \mathbf{1}_{L(\mathcal{K}(\mathbf{k}))}\big)_{\mathbf{i}_{I^c}}\big\Vert_{\mathcal{J}}\prod_{j\in I}\Vert \tilde{Y}_{i_j}^{(j)}\Vert_{L_r}\prod_{j\in I\cap S(\mathbf{k}) } \frac{1}{\sqrt{p_{i_j}}}  \Big)^{(r-2)/r}\nonumber\\
    &\times \Big(e^{-2r} \sum_{\mathbf{i}_{I}}\big\Vert \big(A_d\circ \mathfrak{p}^{(d)}\circ \mathbf{1}_{L(\mathcal{K}(\mathbf{k}))}\big)_{\mathbf{i}_{I^c}}\big\Vert^2_{\mathcal{J}}\prod_{j\in I}\Vert \tilde{Y}_{i_j}^{(j)}\Vert_{L_r}^2\Big)^{1/r}\nonumber\\
    \le& \frac{r-2}{r}\Big(e^{2r/(r-2)}\max_{\mathbf{i}_I}\big\Vert\big(A_d\circ \mathfrak{p}^{(d)}\circ \mathbf{1}_{L(\mathcal{K}(\mathbf{k}))}\big)_{\mathbf{i}_{I^c}}\big\Vert_{\mathcal{J}}\prod_{j\in I}\Vert \tilde{Y}_{i_j}^{(j)}\Vert_{L_r}\prod_{j\in I\cap S(\mathbf{k}) } \frac{1}{\sqrt{p_{i_j}}}  \Big)\nonumber\\
    &+\frac{2}{r} \Big(e^{-2r} \sum_{\mathbf{i}_{I}}\big\Vert \big(A_d\circ \mathfrak{p}^{(d)}\circ \mathbf{1}_{L(\mathcal{K}(\mathbf{k}))}\big)_{\mathbf{i}_{I^c}}\big\Vert^2_{\mathcal{J}}\prod_{j\in I}\Vert \tilde{Y}_{i_j}^{(j)}\Vert_{L_r}^2\Big)^{1/2}\nonumber\\
    \lesssim_{\alpha, d} &\max_{\mathbf{i}_I}\Big(\big\Vert\big(A_d\circ \mathfrak{p}^{(d)}\circ \mathbf{1}_{L(\mathcal{K}(\mathbf{k}))}\big)_{\mathbf{i}_{I^c}}\big\Vert_{\mathcal{J}}\prod_{j\in I}\Vert \tilde{Y}_{i_j}^{(j)}\Vert_{L_r}\prod_{j\in I\cap S(\mathbf{k}) } \frac{1}{\sqrt{p_{i_j}}}\Big)\nonumber\\
    &+\Vert A_d\circ \mathfrak{p}^{(d)}\circ \mathbf{1}_{L(\mathcal{K}(\mathbf{k}))}\Vert_{\mathrm{HS}},
\end{align}
where the last inequality we use the fact that 
\begin{align}
    \prod_{j\in I}\Vert \tilde{Y}_{i_j}^{(j)}\Vert_{L_r}^2\lesssim_{\alpha, d} r^{2d/\alpha},\quad \sum_{\mathbf{i}_{I}}\big\Vert \big(A_d\circ \mathfrak{p}^{(d)}\circ \mathbf{1}_{L(\mathcal{K}(\mathbf{k}))}\big)_{\mathbf{i}_{I^c}}\big\Vert^2_{\mathcal{J}}\le \Vert A_d\circ \mathfrak{p}^{(d)}\circ \mathbf{1}_{L(\mathcal{K}(\mathbf{k}))}\Vert_{\mathrm{HS}}^2.
\end{align}
Note that, when $I=\emptyset$ and $\mathcal{J}=\{[d]\}$, we have 
\begin{align}
    \big\Vert\big(A_d\circ \mathfrak{p}^{(d)}\circ \mathbf{1}_{L(\mathcal{K}(\mathbf{k}))}\big)_{\mathbf{i}_{I^c}}\big\Vert_{\mathcal{J}}\prod_{j\in I}\Vert \tilde{Y}_{i_j}^{(j)}\Vert_{L_r}\prod_{j\in I\cap S(\mathbf{k}) } \frac{1}{\sqrt{p_{i_j}}}= \Vert A_d\circ \mathfrak{p}^{(d)}\circ \mathbf{1}_{L(\mathcal{K}(\mathbf{k}))}\Vert_{\mathrm{HS}}.
\end{align}
Hence, \eqref{Eq_proof_theo1_1} can be further bounded by 
\begin{align}
    & \sum_{d=1}^D\sum_{I\subseteq [d]}\sum_{\mathcal{J}\in \Pi([d]\backslash I)}r^{|\mathcal{J}|/2+|I|/\alpha}\sum_{\nu=1}^d\sum_{\mathbf{k}\in I_{\nu, d}}\max_{\mathbf{i}_I}\Big(\big\Vert\big(A_d\circ \mathfrak{p}^{(d)}\circ \mathbf{1}_{L(\mathcal{K}(\mathbf{k}))}\big)_{\mathbf{i}_{I^c}}\big\Vert_{\mathcal{J}}\prod_{j\in I\cap S(\mathbf{k}) } \frac{1}{\sqrt{p_{i_j}}}\Big)\nonumber\\
    \le&\sum_{d=1}^D\sum_{I\subseteq [d]}\sum_{\mathcal{J}\in \Pi([d]\backslash I)}r^{|\mathcal{J}|/2+|I|/\alpha}\sum_{\nu=1}^d\sum_{\mathbf{k}\in I_{\nu, d}}\max_{\mathbf{i}_I}\Big(\big\Vert\big(A_d\circ \mathfrak{p}^{(d)}\circ \mathbf{1}_{L(\mathcal{K}(\mathbf{k}))}\big)_{\mathbf{i}_{I^c}}\big\Vert_{\mathcal{J}} \frac{1}{\mathfrak{p}^{(d)}_I(\mathbf{i}_I)}\Big)\nonumber\\
    \lesssim_D& \sum_{d=1}^D\sum_{I\subseteq [d]}\sum_{\mathcal{J}\in \Pi([d]\backslash I)}r^{|\mathcal{J}|/2+|I|/\alpha}\max_{\mathbf{i}_I}\Big(\big\Vert\big(A_d\circ \mathfrak{p}^{(d)}\big)_{\mathbf{i}_{I^c}}\big\Vert_{\mathcal{J}} \frac{1}{\mathfrak{p}^{(d)}_I(\mathbf{i}_I)}\Big),
\end{align}
where $\mathfrak{p}^{(d)}_{I}$ was defined in Theorem \ref{Theo_main1} and the second inequality is due to Lemma \ref{Lem_matrixnorm_comparasion}.

To finish the proof, it remains to replace
$\big\Vert\big(A_d\circ \mathfrak{p}^{(d)}\big)_{\mathbf{i}_{I^c}}\big\Vert_{\mathcal{J}}$ by $\big\Vert\big(\mathbb{E}f^{(d)}(X)\circ \mathfrak{p}^{(d)}\big)_{\mathbf{i}_{I^c}}\big\Vert_{\mathcal{J}}.
$
To this end, first fix any multi-index \(\mathbf{i}=(i_1,\cdots, i_d)\in[n]^{d}\) whose distinct values are \(j_{1},\dots ,j_{\nu}\) occurring \(l_{1},\dots ,l_{\nu}\) times, respectively. Then, we have 
\begin{align}\label{Eq_proof_Rd}
    \mathbb{E}\frac{\partial^d{f}}{\partial x_{i_1}\cdots\partial x_{i_d}}(X)=&\sum_{k_1\ge l_1, \cdots, k_{\nu}\ge l_{\nu}}\sum_{m=\nu}^{D}\sum_{\substack{
    k_{\nu +1}, \cdots, k_m>0\\[2pt]    
    k_1+\cdots+k_m\le D}}\sum_{\substack{
    j_{\nu +1}, \cdots, j_m\\[2pt]    
    (j_1,\cdots,j_m)\in [n]^{\underline{m}}}}\nonumber\\
    &\Bigg(\binom{m}{\nu}\nu!c^{(k_1+\cdots+k_{m})}_{(j_1, k_1),\cdots, (j_m, k_m)}\Big(\prod_{\beta=1}^\nu\mathbb{E}X^{k_\beta-l_\beta}_{j_\beta}\frac{k_{\beta}!}{(k_\beta-l_\beta)!}\Big)\Big(\prod_{\beta=\nu+1}^m\mathbb{E}X_{j_{\beta}}^{k_{\beta}}\Big)\Bigg)\nonumber\\
    =&\nu!l_1!\cdots l_\nu!a^{(d)}_{\mathbf{i}}+R_{\mathbf{i}}^{(d)},
\end{align}
where $a_{\mathbf{i}}^{(d)}$ is the $\mathbf{i}$-th entry of $A_{d}$ defined above and $R_{\mathbf{i}}^{(d)}$ corresponds to the set of indices $\mathbf{k}$ such that, $\exists i\in [\nu], k_{i}>l_i$. 

For the case $d=D$, we have $l_1+\cdots+l_\nu=D$. Hence, we have $k_1=l_1, \cdots, k_\nu=l_\nu$ which yields that
\begin{align}\label{Eq_partitial_f}
    \mathbb{E}\frac{\partial^Df}{\partial x_{i_1}\cdots\partial x_{i_d}}(X)=\nu!l_1!\cdots l_\nu!a^{(D)}_{\mathbf{i}}.
\end{align}
Note that $A_D=\sum_{\mathcal{K}\in \Pi_D}A_{D}\circ \mathbf{1}_{L(\mathcal{K})}$. Hence, for any $I\subseteq [d]$ and $\mathcal{J}\in \Pi(I^c)$, we have by the triangle inequality
\begin{align}\label{Eq_main_proof2}
    &\Vert (A_D\circ \mathfrak{p}^{(D)})_{\mathbf{i}_{I^c}}\Vert_{\mathcal{J}}\nonumber\\
    \le& \sum_{\mathcal{K}\in \Pi_D}\Vert (A_{D}\circ \mathfrak{p}^{(D)}\circ\mathbf{1}_{L(\mathcal{K})})_{\mathbf{i}_{I^c}}\Vert_{\mathcal{J}}\le \sum_{\mathcal{K}\in \Pi_D}\Vert (\mathbb{E}f^{(D)}(X)\circ \mathfrak{p}^{(D)}\circ\mathbf{1}_{L(\mathcal{K})})_{\mathbf{i}_{I^c}}\Vert_{\mathcal{J}}\nonumber\\
    \le& \sqrt{2}^{D(D-1)}\cdot|\Pi_D|\cdot\Vert (\mathbb{E}f^{(D)}(X)\circ \mathfrak{p}^{(D)})_{\mathbf{i}_{I^c}}\Vert_{\mathcal{J}}\lesssim_D\Vert (\mathbb{E}f^{(D)}(X)\circ \mathfrak{p}^{(D)})_{\mathbf{i}_{I^c}}\Vert_{\mathcal{J}},
\end{align}
where the third inequality is due to Lemma \ref{Lem_matrixnorm_comparasion} and the second inequality is due to that, for any fixed $\mathcal{K}=\{K_1, \cdots, K_\nu\}\in \Pi_D$, by virtue of \eqref{Eq_partitial_f}, we have 
\begin{align}
    \mathbb{E}f^{(D)}(X)\circ\mathbf{1}_{L(\mathcal{K})}=\nu!|K_1|!\cdots|K_\nu|!\cdot A_D\circ\mathbf{1}_{L(\mathcal{K})}.
\end{align}

To ensure a streamlined presentation, we temporarily assume the validity of the following key estimate, deferring its detailed derivation to Appendix \ref{apendix_supplement}. For every $d=1,\dots ,D-1$, every subset $I\subseteq [d]$, and every partition $\mathcal{J}\in \Pi([d]\setminus I)$, we have
\begin{align}\label{Eq_main_proof3}
    \Vert (R^{(d)}\circ \mathfrak{p}^{(d)})_{\mathbf{i}_{I^c}}\Vert_{\mathcal{J}}\lesssim_D \sum_{k=d+1}^{D}\sum_{\substack{
    \mathcal{K}\in \Pi([k]\backslash I)\\[2pt]    
    |\mathcal{K}|\ge |\mathcal{J}|}}\Vert (A_k\circ \mathfrak{p}^{(k)})_{\mathbf{i}_{I^c}}\Vert_{\mathcal{K}}.
\end{align}

In the specific case where $d=D-1$, invoking \eqref{Eq_proof_Rd}, \eqref{Eq_main_proof2}, and \eqref{Eq_main_proof3} yields
\begin{align}
    &\Vert (A_{D-1}\circ \mathfrak{p}^{(D-1)})_{\mathbf{i}_{I^c}}\Vert_{\mathcal{J}}\nonumber\\
    \lesssim_D& \Vert (\mathbb{E}f^{(D-1)}(X)\circ \mathfrak{p}^{(D-1)})_{\mathbf{i}_{I^c}}\Vert_{\mathcal{J}}+ \Vert (R^{(D-1)}\circ \mathfrak{p}^{(D-1)})_{\mathbf{i}_{I^c}}\Vert_{\mathcal{J}}\nonumber\\
    \lesssim_D&\Vert (\mathbb{E}f^{(D-1)}(X)\circ \mathfrak{p}^{(D-1)})_{\mathbf{i}_{I^c}}\Vert_{\mathcal{J}}+ \sum_{\substack{
    \mathcal{K}\in \Pi([D]\backslash I)\\[2pt]    
    |\mathcal{K}|\ge |\mathcal{J}|}}\Vert (\mathbb{E}f^{(D)}(X)\circ \mathfrak{p}^{(D)})_{\mathbf{i}_{I^c}}\Vert_{\mathcal{K}}.
\end{align}
Consequently, a backward induction argument leads to 
\begin{align}
    &\sum_{d=1}^D\sum_{I\subseteq [d]}\sum_{\mathcal{J}\in \Pi([d]\backslash I)}r^{|\mathcal{J}|/2+|I|/\alpha}\max_{\mathbf{i}_I}\Big(\big\Vert\big(A_d\circ \mathfrak{p}^{(d)}\big)_{\mathbf{i}_{I^c}}\big\Vert_{\mathcal{J}} \frac{1}{\mathfrak{p}^{(d)}_I(\mathbf{i}_I)}\Big)\nonumber\\
    \lesssim_D& \sum_{d=1}^D\sum_{I\subseteq [d]}\sum_{\mathcal{J}\in \Pi([d]\backslash I)}r^{|\mathcal{J}|/2+|I|/\alpha}\max_{\mathbf{i}_I}\Big(\big\Vert\big(\mathbb{E}f^{(d)}(X)\circ \mathfrak{p}^{(d)}\big)_{\mathbf{i}_{I^c}}\big\Vert_{\mathcal{J}} \frac{1}{\mathfrak{p}^{(d)}_I(\mathbf{i}_I)}\Big),
\end{align}
which completes the proof for $r\ge 3$. 

For $1\le r<3$, by monotonicity of $L_r$-norms,
\begin{align}
    \Vert f(X)-\mathbb{E} f(X)\Vert_{L_r}\le \Vert f(X)-\mathbb{E} f(X)\Vert_{L_3}.
\end{align}
Since all powers $r^{\vert \mathcal J\vert/2+\vert I\vert/\alpha}$ are bounded below by constants on $1\le r<3$ and above by constants depending only on $D, \alpha$, the $r=3$ bound implies the desired bound after changing the constants.
\end{proof}

\subsection{Proofs of Theorem \ref{thm:distance} and Corollary \ref{cor:smin}}
\subsubsection{Concentration inequalities for random tensors}\label{subsec:tensor_concentration}

Our proof of the distance theorem (Theorem~\ref{thm:distance}) uses moment estimates for a particular high-order quadratic-type chaos. We first state these estimates in a form that is convenient for later use.

Let \(X^{(1)},\ldots,X^{(d)}\) be independent random vectors in \(\mathbb R^n\). For each \(k\in[d]\) and \(i\in[n]\), assume that
\[
    X_i^{(k)}=\delta_i^{(k)}\xi_i^{(k)},
\]
where \(\delta_i^{(k)}\sim \mathrm{Bernoulli}(p_i)\), \(0<p_i\le 1\), and \(\xi_i^{(k)}\) is mean zero, has variance one, and is \(\alpha\)-sub-exponential with
\[
    \|\xi_i^{(k)}\|_{\Psi_\alpha}\le 1
\]
for some fixed \(0<\alpha\le 1\). All random variables appearing above are assumed to be independent. We write
\[
    X=(X^{(1)},\ldots,X^{(d)})\in\mathbb R^{nd}
\]
for the concatenated random vector.

Let \(A=(A_{i_1,\ldots,i_{2d}})_{(i_1,\ldots,i_{2d})\in[n]^{2d}}\) be a deterministic \(2d\)-tensor, and consider the polynomial
\begin{equation}\label{eq:tensor_poly}
    f(X)
    =
    \sum_{i_1,\ldots,i_{2d}=1}^{n}
    A_{i_1,\ldots,i_d,i_{d+1},\ldots,i_{2d}}
    \prod_{k=1}^{d} X^{(k)}_{i_k}X^{(k)}_{i_{d+k}} .
\end{equation}
This polynomial is homogeneous of degree two in each block \(X^{(k)}\), and has total degree \(2d\). In the proof of Theorem~\ref{thm:distance}, \(A\) will be chosen so that \(f(X)\) represents the squared norm of a projected random tensor. Thus the moment estimates for \eqref{eq:tensor_poly} provide the input needed for the distance theorem.

To apply the general polynomial concentration result, Theorem~\ref{Theo_main1}, the main technical point is to estimate the expected derivative tensors
\[
    \mathbb E f^{(m)}(X), \qquad 1\le m\le 2d,
\]
in the relevant tensor partition norms. The block-homogeneous structure of \eqref{eq:tensor_poly} allows these norms to be written in more concrete terms, using matrix operator norms and certain partial traces of the tensor \(A\). We therefore introduce the notation and weights that will be used in the statement and proof of the main result of this subsection.

We identify the coordinate $X_i^{(k)}$ with the block-coordinate pair $(k,i)$.
Equivalently, one may identify $(k,i)$ with $(k-1)n+i\in[nd]$. For a multi-index
\[
    \gamma=((k_1,i_1),\ldots,(k_m,i_m))\in([d]\times[n])^m,
\]
define the sparse weight
\[
    \mathfrak{P}^{(m)}_\gamma
    =
    \left(
    \prod_{(k,i)\in \operatorname{supp}(\gamma)} p_i
    \right)^{1/2},
\]
where $\operatorname{supp}(\gamma)$ denotes the set of distinct block-coordinate
pairs appearing in $\gamma$. For $S\subseteq[m]$, the marginal weight is
\[
    \mathfrak{P}^{(m)}_S(\gamma_S)=
    \left(
    \prod_{(k,i)\in \operatorname{supp}(\gamma_S)} p_i
    \right)^{1/2},
\]
where $\gamma_S=((k_j, i_j), j\in S).$ This is the same weight as in Theorem~1 after the identification
$(k,i)\leftrightarrow (k-1)n+i$.

\paragraph{Reduced weights for the free variables.}
For $l\in[d]$ and $j=(j_1,\ldots,j_{2l})\in[n]^{2l}$, define the reduced
block-coordinate map
\[
    \vartheta_l(s,j_s)
    :=
    \begin{cases}
        (s,j_s), & 1\le s\le l,\\
        (s-l,j_s), & l+1\le s\le 2l.
    \end{cases}
\]
For $I_2\subseteq[2l]$, define
\begin{equation}\label{eq:tilde_P_marginal}
    \widetilde{\mathfrak{P}}_{I_2}^{(2l)}(j_{I_2})
    :=
    \left(
    \prod_{(m,a)\in
    \{\vartheta_l(s,j_s):\,s\in I_2\}}
    p_a
    \right)^{1/2}.
\end{equation}
The product is taken over distinct pairs $(m,a)$. In particular, if two equal
coordinates occur in the same reduced block, they contribute only one factor
$\sqrt p$; if the same coordinate value occurs in two different reduced blocks,
it contributes twice, since the underlying random variables belong to different
independent copies.

We write $\widetilde{\mathfrak{P}}^{(2l)}(j):=\widetilde{\mathfrak{P}}_{[2l]}^{(2l)}(j)$.
Equivalently, for $j=(a,b)=(a_1,\ldots,a_l,b_1,\ldots,b_l)$,
\begin{equation}\label{eq:tilde_P_full}
    \widetilde{\mathfrak{P}}^{(2l)}(a,b)
    =
    \prod_{m=1}^l q(a_m,b_m),
    \qquad
    q(a_m,b_m)
    :=
    \begin{cases}
        \sqrt{p_{a_m}p_{b_m}}, & a_m\ne b_m,\\[2mm]
        \sqrt{p_{a_m}}, & a_m=b_m.
    \end{cases}
\end{equation}

\paragraph{Weighted partial traces.}
Fix $I=\{k_1,\ldots,k_l\}\subseteq[d]$, where $k_1<\cdots<k_l$.
For $a=(a_1,\ldots,a_l)\in[n]^l$, $b=(b_1,\ldots,b_l)\in[n]^l$,
$u=(u_t)_{t\in[d]\setminus I}\in[n]^{d-l}$, and
$\varepsilon=(\varepsilon_1,\ldots,\varepsilon_l)\in\{0,1\}^l$, define
\[
    \Theta_{I,\varepsilon}(a,b,u)
    =
    (i_1,\ldots,i_{2d})\in[n]^{2d}
\]
as follows. For each unselected block $t\in[d]\setminus I$, set
\[
    i_t=i_{d+t}=u_t.
\]
For each selected block $t=k_m\in I$, set
\[
    (i_{k_m},i_{d+k_m})
    =
    \begin{cases}
        (a_m,b_m), & \varepsilon_m=0,\\
        (b_m,a_m), & \varepsilon_m=1.
    \end{cases}
\]
Define the $p$-weighted partial trace
\begin{equation}\label{eq:T_I_p_A}
    T_I^p A(a,b)
    :=
    \sum_{u\in[n]^{d-l}}
    \left(\prod_{t\in[d]\setminus I}p_{u_t}\right)
    \sum_{\varepsilon\in\{0,1\}^l}
    A_{\Theta_{I,\varepsilon}(a,b,u)}.
\end{equation}
Finally, define the reduced tensor
\begin{equation}\label{eq:corrected_reduced_tensor}
    \mathcal A_I^{(2l,p)}(a,b)
    :=
    \widetilde{\mathfrak{P}}^{(2l)}(a,b)\,T_I^p A(a,b).
\end{equation}
The factor $\prod_{t\notin I}p_{u_t}$ comes from the contraction
\[
    \mathbb E X_{i_t}^{(t)}X_{i_{d+t}}^{(t)}
    =
    p_{i_t}\mathbf 1_{\{i_t=i_{d+t}\}},
\]
whereas $\widetilde{\mathfrak{P}}^{(2l)}(a,b)$ is the derivative weight from Theorem~\ref{Theo_main1}
for the free variables.

\begin{theorem}
\label{thm:tensor_poly_concentration}
Let $f$ be the polynomial defined in \eqref{eq:tensor_poly}. Then, for every
$r\ge 2$,
\begin{equation}\label{eq:thm3_corrected}
\begin{aligned}
    \|f(X)-\mathbb Ef(X)\|_{L_r}
    \lesssim_{\alpha,d}
    \sum_{l=1}^d
    \sum_{I_2\subseteq[2l]}
    \sum_{\mathcal J\in\Pi([2l]\setminus I_2)}
    r^{|\mathcal J|/2+|I_2|/\alpha}
    \max_{j_{I_2}}
    \sum_{\substack{I_1\subseteq[d]\\ |I_1|=l}}
    \frac{
    \left\|
    \left(\mathcal A_{I_1}^{(2l,p)}\right)_{j_{I_{2}^{c}}}
    \right\|_{\mathcal J}
    }{
    \widetilde{\mathfrak{P}}_{I_2}^{(2l)}(j_{I_2})
    } .
\end{aligned}
\end{equation}
 We use the conventions
$\widetilde{\mathfrak{P}}_{\varnothing}^{(2l)}=1$ and
$\Pi(\varnothing)=\{\varnothing\}$.
\end{theorem}

\begin{proof}
Applying Theorem~\ref{Theo_main1} to this polynomial, whose
degree is at most $2d$, gives
\begin{equation}\label{eq:apply_thm1_to_tensor_poly}
\begin{aligned}
    \|f(X)-\mathbb Ef(X)\|_{L_r}
    \lesssim_{\alpha,d}
    \sum_{m=1}^{2d}
    \sum_{S\subseteq[m]}
    \sum_{\mathcal J\in\Pi([m]\setminus S)}
    r^{|\mathcal J|/2+|S|/\alpha}
    \max_{\gamma_S}
    \frac{
    \left\|
    \left(\mathbb Ef^{(m)}(X)\circ \mathfrak P^{(m)}\right)_{\gamma_{S^{c}}}
    \right\|_{\mathcal J}
    }{
    \mathfrak P_S^{(m)}(\gamma_S)
    } .
\end{aligned}
\end{equation}
Here $\gamma=((k_1,i_1),\ldots,(k_m,i_m))$ is a derivative multi-index in
$([d]\times[n])^m$.

We next identify the nonzero entries of $\mathbb Ef^{(m)}(X)$. For a derivative
multi-index $\gamma$, let
\[
    N_t(\gamma)
    :=
    \#\{s\in[m]: k_s=t\},
    \qquad t\in[d].
\]
Since the polynomial \eqref{eq:tensor_poly} has degree at most two in each
block $X^{(t)}$, we have
\[
    \mathbb E\partial_\gamma f(X)=0
\]
unless
\[
    N_t(\gamma)\in\{0,2\}
    \qquad\text{for every }t\in[d].
\]
Indeed, if $N_t(\gamma)\ge 3$, the derivative with respect to the $t$-th block
vanishes. If $N_t(\gamma)=1$, then after differentiation one centered factor
from the $t$-th block remains, and its expectation is zero.

Consequently only even orders $m=2l$ contribute. Let
\[
    I=\{k_1,\ldots,k_l\}\subseteq[d],
    \qquad k_1<\cdots<k_l,
\]
be the set of active blocks, namely the blocks differentiated exactly twice.
For $a=(a_1,\ldots,a_l)$ and $b=(b_1,\ldots,b_l)$, consider the canonical
derivative
\[
    \partial_{I,a,b}
    :=
    \partial_{z_{k_1,a_1}}\cdots\partial_{z_{k_l,a_l}}
    \partial_{z_{k_1,b_1}}\cdots\partial_{z_{k_l,b_l}} .
\]
For each selected block $k_m$,
\[
\begin{aligned}
    &\partial_{z_{k_m,a_m}}\partial_{z_{k_m,b_m}}
    \left(
    z_{k_m,i_{k_m}}z_{k_m,i_{d+k_m}}
    \right)  \\
    &\qquad =
    \mathbf 1_{\{(i_{k_m},i_{d+k_m})=(a_m,b_m)\}}
    +
    \mathbf 1_{\{(i_{k_m},i_{d+k_m})=(b_m,a_m)\}} .
\end{aligned}
\]
This identity also gives the correct factor $2$ when $a_m=b_m$.
For each unselected block $t\notin I$,
\[
    \mathbb E
    X_{i_t}^{(t)}X_{i_{d+t}}^{(t)}
    =
    p_{i_t}\mathbf 1_{\{i_t=i_{d+t}\}} .
\]
Therefore
\begin{equation}\label{eq:derivative_expectation_formula}
\begin{aligned}
    \mathbb E\,\partial_{I,a,b} f(X)
    &=
    \sum_{u\in[n]^{d-l}}
    \left(\prod_{t\in[d]\setminus I}p_{u_t}\right)
    \sum_{\varepsilon\in\{0,1\}^l}
    A_{\Theta_{I,\varepsilon}(a,b,u)}        \\
    &= T_I^pA(a,b).
\end{aligned}
\end{equation}
Multiplying \eqref{eq:derivative_expectation_formula} by the derivative weight
from Theorem~1 yields
\begin{equation}\label{eq:weighted_derivative_reduced_tensor}
    \left(\mathbb Ef^{(2l)}(X)\circ \mathfrak P^{(2l)}\right)_{I,a,b}
    =
    \widetilde{\mathfrak{P}}^{(2l)}(a,b)\,T_I^pA(a,b)
    =
    \mathcal A_I^{(2l,p)}(a,b),
\end{equation}
for the canonical ordering of the derivative coordinates.

It remains only to account for the fact that the tensor
$\mathbb Ef^{(2l)}(X)\circ \mathfrak P^{(2l)}$ is indexed by all possible orderings of
the $2l$ derivative coordinates. Every nonzero entry is obtained from one of
the canonical tensors $\mathcal A_I^{(2l,p)}$, $|I|=l$, by a coordinate
permutation and by the corresponding embedding into the appropriate block
coordinates. The number of such permutations and embeddings depends only on
$d$.

Partition norms are invariant under coordinate relabeling, up to relabeling
the underlying partition. Moreover, under the same relabeling the denominator
$\mathfrak P_S^{(2l)}(\gamma_S)$ becomes precisely the reduced marginal weight
$\widetilde{\mathfrak{P}}_{\pi(S)}^{(2l)}(j_{\pi(S)})$ for the corresponding subset of
coordinates. Since the right-hand side of \eqref{eq:thm3_corrected} sums over
all subsets $I_2\subseteq[2l]$ and all partitions
$\mathcal J\in\Pi([2l]\setminus I_2)$, these finitely many relabelings
are absorbed into the constant depending only on $d$. Hence, for each
$l=1,\ldots,d$, the total contribution of the $m=2l$ term in
\eqref{eq:apply_thm1_to_tensor_poly} is bounded by
\[
\begin{aligned}
    C_d
    \sum_{I_2\subseteq[2l]}
    \sum_{\mathcal J\in\Pi([2l]\setminus I_2)}
    r^{|\mathcal J|/2+|I_2|/\alpha}
    \max_{j_{I_2}}
    \sum_{\substack{I_1\subseteq[d]\\ |I_1|=l}}
    \frac{
    \left\|
    \left(\mathcal A_{I_1}^{(2l,p)}\right)_{j_{I_{2}^{c}}}
    \right\|_{\mathcal J}
    }{
    \widetilde{\mathfrak{P}}_{I_2}^{(2l)}(j_{I_2})
    } .
\end{aligned}
\]
Summing over $l=1,\ldots,d$ proves \eqref{eq:thm3_corrected}.
\end{proof}

While Theorem \ref{thm:tensor_poly_concentration} gives a precise moment bound, the sum over all partitions is cumbersome in applications. We therefore specialize to the homogeneous setting, where the selection probabilities are uniform. Bounding the partition norms by tensor operator norms then yields the following corollary. This bound separates the sparse Gaussian-type fluctuation from the heavy-tailed fractional-exponential decay, and it is the main tool in the proof of our subspace distance theorem.

% Corrected Corollary 1 and Appendix B
% This file is intended to replace the old Corollary 1 and the old Appendix B.
% It assumes the notation of the corrected Theorem 3, in particular
% the corrected tensors \mathcal A_I^{(2l,p)} and the weights \widetilde{\mathfrak{P}}_{I_2}^{(2l)}.

\begin{coro}
\label{cor:cor1_corrected_refined}
Assume the setting of Theorem~\ref{thm:tensor_poly_concentration}. Suppose in addition that
\[
    p_1=\cdots=p_n=p\in(0,1].
\]
Define
\begin{equation}\label{eq:Gamma_npd_def}
    \Gamma_{n,p,d}
    :=
    p^{d/2}\max_{0\le s\le d-1}(np)^{s/2}.
\end{equation}
Then, for every $r\ge2$,
\begin{equation}\label{eq:cor1_corrected_refined}
\begin{aligned}
    \|f(X)-\mathbb Ef(X)\|_{L_r}
    \lesssim_{\alpha,d}
    &\ r^{1/2}\Gamma_{n,p,d}\,
    \|A\|_{[2d]}                                      \\
    &\quad+
    \sum_{l=1}^dp^{d-l}
    \sum_{k=2}^{2l}
    r^{k/\alpha}n^{d-k/2}
    \|A\|_{[d],[2d]\setminus[d]}.
\end{aligned}
\end{equation}
\end{coro}

\begin{proof}
We apply Theorem~\ref{thm:tensor_poly_concentration}. Fix
$l\in[d]$, $I_1\subseteq[d]$ with $|I_1|=l$, a subset
$I_2\subseteq[2l]$, and a partition
$\mathcal J\in\Pi([2l]\setminus I_2)$. Put
\[
    q:=|I_2|,
    \qquad
    \kappa:=|\mathcal J|.
\]
In the homogeneous case, the reduced tensor admits the decomposition
\begin{equation}\label{eq:cor1_homogeneous_decomposition}
    \mathcal A_{I_1}^{(2l,p)}
    =
    p^{d-l}\,\widetilde{\mathfrak{P}}^{(2l)}
    \circ
    \sum_{\varepsilon\in\{0,1\}^l}T_{I_1,\varepsilon}A,
\end{equation}
where $T_{I_1,\varepsilon}A$ is the unweighted partial trace defined in
Appendix~\ref{app:cor1_appendixB}. We split the contribution in
Theorem~\ref{thm:tensor_poly_concentration} according to whether
$\kappa+q=1$ or $\kappa+q\ge2$.

First consider the terms with $\kappa+q=1$. Since $\mathcal J$ is a partition
of $[2l]\setminus I_2$, this can occur only when $I_2=\varnothing$ and
$\mathcal J=\{[2l]\}$. Hence the norm is the Hilbert--Schmidt norm. By
\eqref{eq:cor1_homogeneous_decomposition}, the pointwise bound
$\widetilde{\mathfrak{P}}^{(2l)}(j)\le p^{l/2}$, and
\eqref{eq:B_partial_trace_hs_estimate},
\begin{align}
    \left\|\mathcal A_{I_1}^{(2l,p)}\right\|_{[2l]}
    &\le
    p^{d-l}
    \sum_{\varepsilon\in\{0,1\}^l}
    \left\|
    \widetilde{\mathfrak{P}}^{(2l)}\circ T_{I_1,\varepsilon}A
    \right\|_{[2l]}                                      \notag\\
    &\le
    p^{d-l/2}
    \sum_{\varepsilon\in\{0,1\}^l}
    \left\|T_{I_1,\varepsilon}A\right\|_{[2l]}             \notag\\
    &\lesssim_d
    p^{d-l/2}n^{(d-l)/2}\|A\|_{[2d]} .
    \label{eq:cor1_HS_contribution_bound}
\end{align}
Therefore, after summing over the finitely many choices of $l$ and $I_1$, the
total contribution of the terms with $\kappa+q=1$ is bounded by
\begin{align}
    C_{\alpha,d}r^{1/2}
    \max_{1\le l\le d}
    p^{d-l/2}n^{(d-l)/2}\|A\|_{[2d]}
    &=
    C_{\alpha,d}r^{1/2}
    p^{d/2}\max_{0\le s\le d-1}(np)^{s/2}\|A\|_{[2d]}       \notag\\
    &=
    C_{\alpha,d}r^{1/2}\Gamma_{n,p,d}\|A\|_{[2d]} .
    \label{eq:cor1_first_part}
\end{align}

It remains to handle the terms with $\kappa+q\ge2$. We first remove the free
weight $\widetilde{\mathfrak{P}}^{(2l)}$. If $I_2=\varnothing$, then
\eqref{eq:B_free_weight_multiplier_empty} and $p\le1$ give
\begin{align}
    \left\|\mathcal A_{I_1}^{(2l,p)}\right\|_{\mathcal J}
    &\le
    p^{d-l}
    \sum_{\varepsilon\in\{0,1\}^l}
    \left\|
    \widetilde{\mathfrak{P}}^{(2l)}\circ T_{I_1,\varepsilon}A
    \right\|_{\mathcal J}                                  \notag\\
    &\lesssim_d
    p^{d-l/2}
    \sum_{\varepsilon\in\{0,1\}^l}
    \left\|T_{I_1,\varepsilon}A\right\|_{\mathcal J}.
    \label{eq:cor1_remove_weight_empty}
\end{align}
If $I_2\ne\varnothing$, then \eqref{eq:B_free_weight_multiplier_sliced} gives,
uniformly in $j_{I_2}$,
\begin{align}
    \frac{
    \left\|
    \left(\mathcal A_{I_1}^{(2l,p)}\right)_{j_{I_{2}^{c}}}
    \right\|_{\mathcal J}
    }{
    \widetilde{\mathfrak{P}}_{I_2}^{(2l)}(j_{I_2})
    }
    &\le
    p^{d-l}
    \sum_{\varepsilon\in\{0,1\}^l}
    \frac{
    \left\|
    \left(\widetilde{\mathfrak{P}}^{(2l)}\circ T_{I_1,\varepsilon}A\right)_{j_{I_{2}^{c}}}
    \right\|_{\mathcal J}
    }{
    \widetilde{\mathfrak{P}}_{I_2}^{(2l)}(j_{I_2})
    }                                                        \notag\\
    &\lesssim_d
    p^{d-l}
    \sum_{\varepsilon\in\{0,1\}^l}
    \left\|
    \left(T_{I_1,\varepsilon}A\right)_{j_{I_{2}^{c}}}
    \right\|_{\mathcal J}.
    \label{eq:cor1_remove_weight_nonempty}
\end{align}
Combining \eqref{eq:cor1_remove_weight_empty} and
\eqref{eq:cor1_remove_weight_nonempty}, and then applying
\eqref{eq:B_partial_trace_op_estimate}, we obtain, for all terms with
$\kappa+q\ge2$,
\begin{equation}\label{eq:cor1_operator_reduced_bound}
    \frac{
    \left\|
    \left(\mathcal A_{I_1}^{(2l,p)}\right)_{j_{I_{2}^{c}}}
    \right\|_{\mathcal J}
    }{
    \widetilde{\mathfrak{P}}_{I_2}^{(2l)}(j_{I_2})
    }
    \lesssim_d p^{d-l}
    n^{d-(\kappa+q)/2}
    \|A\|_{[d],[2d]\setminus[d]} .
\end{equation}
Here, when $I_2=\varnothing$, the denominator is interpreted as one.

The corresponding factor in Theorem~\ref{thm:tensor_poly_concentration} is
$r^{\kappa/2+q/\alpha}$. Since $0<\alpha\le1$ and $r\ge2$,
\begin{equation}\label{eq:cor1_r_power_bound}
    r^{\kappa/2+q/\alpha}
    \le
    r^{(\kappa+q)/\alpha}.
\end{equation}
Set $k:=\kappa+q$. Since we are in the case $k\ge2$ and always have
$k\le2d$, summing over all choices of $l$, $I_1$, $I_2$, and $\mathcal J$
only changes the implicit constant by a factor depending on $d$. Hence the
total contribution of the terms with $\kappa+q\ge2$ is bounded by
\begin{equation}\label{eq:cor1_second_part}
    C_{\alpha,d}\sum_{l=1}^d p^{d-l}
    \sum_{k=2}^{2l}
    r^{k/\alpha}n^{d-k/2}
    \|A\|_{[d],[2d]\setminus[d]} .
\end{equation}
Combining \eqref{eq:cor1_first_part} and \eqref{eq:cor1_second_part} proves
\eqref{eq:cor1_corrected_refined}.
\end{proof}

\subsubsection{Distance to subspaces for random tensors}

\begin{proof}[Proof of Theorem \ref{thm:distance}]
Assume $B\neq0$ and set
\[
    M:=B^\top B.
\]
We view $M$ as a $2d$-tensor by writing
\[
    A_{i_1,\ldots,i_d,i_{d+1},\ldots,i_{2d}}
    :=
    M_{(i_1,\ldots,i_d),(i_{d+1},\ldots,i_{2d})}.
\]
Then
\[
    X^\top M X
    =
    \sum_{i_1,\ldots,i_{2d}=1}^n
    A_{i_1,\ldots,i_{2d}}
    \prod_{k=1}^d X_{i_k}^{(k)}X_{i_{d+k}}^{(k)}.
\]

We first record the quadratic-form consequence of Corollary~\ref{cor:cor1_corrected_refined}. For every
$r\ge2$,
\begin{equation}\label{eq:corrected_quadratic_moment}
\begin{aligned}
    \left\|
    X^\top M X-\mathbb E X^\top M X
    \right\|_{L_r}
    \lesssim_{\alpha,d}
    r^{1/2}\Gamma_{n,p,d}\|M\|_{\mathrm{HS}}
    +
    \sum_{k=2}^{2d} r^{k/\alpha}n^{d-k/2}\|M\|_{\mathrm{op}}.
\end{aligned}
\end{equation}

Since the coordinates of the simple tensor are isotropic up to the factor $p^d$, we have
\begin{equation}\label{eq:expectation_quadratic_form}
    \mathbb E X^\top M X
    =p^d\operatorname{tr}(M)
    =p^d\|B\|_{\mathrm{HS}}^2.
\end{equation}
Put
\[
    Y:=X^\top M X-\mathbb E X^\top M X,
    \qquad
    b:=p^{d/2}\|B\|_{\mathrm{HS}}.
\]
For non-negative real numbers $a,b$ with $b>0$,
\begin{equation}\label{eq:sqrt_algebraic_ineq}
    |a-b|
    \le
    \min\left\{
    \frac{|a^2-b^2|}{b},\sqrt{|a^2-b^2|}
    \right\}.
\end{equation}
Applying this with $a=\|BX\|_2$ and using \eqref{eq:expectation_quadratic_form}, we obtain
\[
    \left\|
    \|BX\|_2-p^{d/2}\|B\|_{\mathrm{HS}}
    \right\|_{L_r}
    \le
    \min\left\{
    \frac{\|Y\|_{L_r}}{b},
    \left\|\sqrt{|Y|}\right\|_{L_r}
    \right\}.
\]
Since $r\ge2$,
\[
    \left\|\sqrt{|Y|}\right\|_{L_r}
    =\|Y\|_{L_{r/2}}^{1/2}
    \le
    \|Y\|_{L_r}^{1/2}.
\]
Let
\[
    g(x):=\min\{x/b,\sqrt x\},\qquad x\ge0.
\]
The function $g$ is non-negative, increasing, concave, and satisfies $g(0)=0$, hence it is
subadditive. Applying $g$ to the right-hand side of
\eqref{eq:corrected_quadratic_moment}, using $g(x)\le x/b$ for the Hilbert--Schmidt term and
$g(x)\le\sqrt{x}$ for the remaining terms, we get
\begin{equation}\label{eq:distance_moment_corrected}
\begin{aligned}
    \left\|
    \|BX\|_2-p^{d/2}\|B\|_{\mathrm{HS}}
    \right\|_{L_r}
    \lesssim_{\alpha,d}
    &\ r^{1/2}
    \frac{\Gamma_{n,p,d}}{p^{d/2}}
    \frac{\|B^\top B\|_{\mathrm{HS}}}{\|B\|_{\mathrm{HS}}} \\
    &\quad+
    \sum_{k=2}^{2d}
    r^{k/(2\alpha)}n^{d/2-k/4}
    \sqrt{\|B^\top B\|_{\mathrm{op}}}.
\end{aligned}
\end{equation}
Using
\[
    \|B^\top B\|_{\mathrm{HS}}
    \le
    \|B\|_{\mathrm{op}}\|B\|_{\mathrm{HS}},
    \qquad
    \|B^\top B\|_{\mathrm{op}}=\|B\|_{\mathrm{op}}^2,
\]
 this becomes
\begin{equation}\label{eq:distance_moment_corrected_simplified}
\begin{aligned}
    \left\|
    \|BX\|_2-p^{d/2}\|B\|_{\mathrm{HS}}
    \right\|_{L_r}
    \lesssim_{\alpha,d}
    \|B\|_{\mathrm{op}}
    \left(
    r^{1/2}\Lambda_{n,p,d}
    +
    \sum_{k=2}^{2d}
    r^{k/(2\alpha)}n^{d/2-k/4}
    \right).
\end{aligned}
\end{equation}

By the standard moment-to-tail conversion lemma, see, e.g., \cite[Lemma~2.1]{DSUW2025JFA}, \eqref{eq:distance_moment_corrected_simplified} implies that, for every $t>0$,
\begin{equation}\label{eq:tail_before_phase_simplification}
\begin{aligned}
    \mathbb P\left\{
    \left|
    \|BX\|_2-p^{d/2}\|B\|_{\mathrm{HS}}
    \right|>t
    \right\}
    \le
    2\exp\Bigg(-c(\alpha,d)
    \min\Bigg\{
    &\frac{t^2}{\Lambda_{n,p,d}^2\|B\|_{\mathrm{op}}^2},\\
    &\min_{2\le k\le2d}
    \left(
    \frac{t}{n^{d/2-k/4}\|B\|_{\mathrm{op}}}
    \right)^{2\alpha/k}
    \Bigg\}
    \Bigg).
\end{aligned}
\end{equation}
It remains to simplify the second minimum. Let
\[
    x:=\frac{t}{n^{d/2}\|B\|_{\mathrm{op}}}.
\]
For $2\le k\le2d$,
\[
    \left(
    \frac{t}{n^{d/2-k/4}\|B\|_{\mathrm{op}}}
    \right)^{2\alpha/k}
    =
    n^{\alpha/2}x^{2\alpha/k}.
\]
If $x\le1$, the minimum over $2\le k\le2d$ is attained at $k=2$, giving
\[
    n^{\alpha/2}x^\alpha
    =
    \left(
    \frac{t}{n^{(d-1)/2}\|B\|_{\mathrm{op}}}
    \right)^\alpha.
\]
If $x>1$, the minimum is attained at $k=2d$, giving
\[
    n^{\alpha/2}x^{\alpha/d}
    =
    \left(
    \frac{t}{\|B\|_{\mathrm{op}}}
    \right)^{\alpha/d}.
\]
Thus the second minimum in \eqref{eq:tail_before_phase_simplification} is bounded below by
\[
    \min\left\{
    \left(\frac{t}{n^{(d-1)/2}\|B\|_{\mathrm{op}}}\right)^\alpha,
    \left(\frac{t}{\|B\|_{\mathrm{op}}}\right)^{\alpha/d}
    \right\}.
\]
\end{proof}

\subsubsection{The Singularity of simple tensor matrices}

\begin{proof}[Proof of Corollary~\ref{cor:smin}]
Let \(N=n^d\).
For each \(j\in[m]\), set
\[
    H_j:=\operatorname{span}\{\mathbf X_i:i\ne j\},
    \qquad
    Q_j:=P_{H_j^\perp},
    \qquad
    k_j:=\operatorname{rank}(Q_j).
\]
Since \(m\le (1-\varepsilon)N\), we have
\[
    k_j=\operatorname{codim}(H_j)\ge N-m+1\ge \varepsilon N .
\]
Put
\[
    \eta:=p^d\varepsilon n .
\]
By the assumption on \(\varepsilon\), after choosing \(c_1=c_1(\alpha,d)\) sufficiently large,
we may assume that \(\eta\ge 1\). Moreover, since \(\varepsilon<1\), the range of parameters is
non-empty only if \(np^d\ge c_1(d\log n)^{2/\alpha}\). Hence, increasing \(c_1\) if necessary,
\(np\ge np^d\ge 1\), and therefore
\[
    \Lambda_{n,p,d}=\max_{0\le s\le d-1}(np)^{s/2}
    =(np)^{(d-1)/2}.
\]

Conditioning on all columns except \(\mathbf X_j\), the projection \(Q_j\) is fixed and \(\mathbf X_j\) is
independent of \(Q_j\). Applying Theorem~\ref{thm:distance} with \(B=Q_j\) and using
\[
    \|Q_j\|_{\mathrm{op}}=1,
    \qquad
    \|Q_j\|_{\mathrm{HS}}=\sqrt{k_j},
\]
we get
\[
\begin{aligned}
&\mathbb P\left\{
    \operatorname{dist}(\mathbf X_j,H_j)
    \le \frac12 p^{d/2}\sqrt{\varepsilon N}
    \,\middle|\, H_j
\right\}                                                     \\
&\qquad\le
\mathbb P\left\{
    \|Q_j\mathbf X_j\|_2
    \le \frac12 p^{d/2}\sqrt{k_j}
    \,\middle|\, H_j
\right\}                                                     \\
&\qquad\le
2\exp\left[
    -c(\alpha,d)
    \min\left\{
        \frac{p^d k_j}{(np)^{d-1}},
        \left(
            \frac{p^{d/2}\sqrt{k_j}}{n^{(d-1)/2}}
        \right)^\alpha,
        \left(
            p^{d/2}\sqrt{k_j}
        \right)^{\alpha/d}
    \right\}
\right].
\end{aligned}
\]
Since \(k_j\ge \varepsilon n^d\), \(0<p\le1\), \(0<\varepsilon<1\), and \(\eta\ge1\), the three
terms in the minimum are bounded from below, up to constants depending only on \(\alpha,d\), by
\[
    p\varepsilon n,
    \qquad
    (p^d\varepsilon n)^{\alpha/2},
    \qquad
    p^{\alpha/2}\varepsilon^{\alpha/(2d)}n^{\alpha/2}.
\]
Each of these is at least \(c_{\alpha,d}\eta^{\alpha/2}\). Indeed,
\[
    p\varepsilon n\ge p^d\varepsilon n=\eta\ge \eta^{\alpha/2},
\]
and
\[
    p^{\alpha/2}\varepsilon^{\alpha/(2d)}n^{\alpha/2}
    \ge
    p^{d\alpha/2}\varepsilon^{\alpha/2}n^{\alpha/2}
    =
    \eta^{\alpha/2}.
\]
Therefore,
\[
    \mathbb P\left\{
        \operatorname{dist}(\mathbf X_j,H_j)
        \le \frac12 p^{d/2}\sqrt{\varepsilon N}
    \right\}
    \le
    2\exp\left(-c_{\alpha,d}\eta^{\alpha/2}\right).
\]
Taking a union bound over \(j=1,\ldots,m\), and using \(m\le N=n^d\), gives
\[
\begin{aligned}
    \mathbb P\left\{
        \exists j\in[m]:
        \operatorname{dist}(\mathbf X_j,H_j)
        \le \frac12 p^{d/2}\sqrt{\varepsilon N}
    \right\}
    &\le
    2n^d\exp\left(-c_{\alpha,d}\eta^{\alpha/2}\right)        \\
    &\le
    2\exp\left(-c_2(\alpha,d)\eta^{\alpha/2}\right),
\end{aligned}
\]
where the last inequality follows from
\[
    \eta^{\alpha/2}
    =
    (p^d\varepsilon n)^{\alpha/2}
    \ge
    c_1^{\alpha/2}d\log n
\]
and from choosing \(c_1\) sufficiently large.

On the complementary event,
\[
    \operatorname{dist}(\mathbf X_j,H_j)
    >
    \frac12 p^{d/2}\sqrt{\varepsilon N},
    \qquad j=1,\ldots,m.
\]
Let $s_{\ell}(\mathbf X)$ denote the singular values of $\mathbf X$. By the negative second moment identity,
\[
    \sum_{\ell=1}^m s_\ell(\mathbf X)^{-2}
    =
    \sum_{j=1}^m \operatorname{dist}(X_j,H_j)^{-2}.
\]
Thus
\[
\begin{aligned}
    \sigma_{\min}(\mathbf X)^{-2}
    \le
    \sum_{\ell=1}^m s_\ell(\mathbf X)^{-2}               =
    \sum_{j=1}^m \operatorname{dist}(\mathbf X_j,H_j)^{-2}        
    <
    \frac{4m}{p^d\varepsilon N}
    \le
    \frac{4(1-\varepsilon)}{p^d\varepsilon}.
\end{aligned}
\]
Equivalently,
\[
    \sigma_{\min}(\mathbf X)
    >
    \frac12 p^{d/2}\sqrt{\frac{\varepsilon}{1-\varepsilon}}.
\]
Since \(\eta=p^d\varepsilon n\), the desired estimate follows:
\[
    \mathbb P\left\{
        \sigma_{\min}(\mathbf X)
        >
        \frac12 p^{d/2}\sqrt{\frac{\varepsilon}{1-\varepsilon}}
    \right\}
    \ge
    1-2\exp\left(-c_2(\alpha,d)(p^d\varepsilon n)^{\alpha/2}\right).
\]
\end{proof}

\appendix
\section{Supplementary Proof of Theorem \ref{Theo_main1}} \label{apendix_supplement}

In this appendix we prove the estimate used in the proof of Theorem 1:
for every \(d=1,\ldots,D-1\), every \(I\subseteq[d]\), every fixed \(\mathbf{i}_I\), and every
\(\mathcal{J}\in \Pi([d]\setminus I)\),
\begin{equation}\label{A.1}
\left\|
(R^{(d)}\circ \mathfrak{p}^{(d)})_{\mathbf{i}_{I^{c}}}
\right\|_\mathcal{J}
\lesssim_{\alpha,D}
\sum_{q=d+1}^{D}
\sum_{\substack{\mathcal{K}\in \Pi([q]\setminus I)\\ |\mathcal{K}|\ge |\mathcal{J}|}}
\left\|
(A_q\circ \mathfrak{p}^{(q)})_{\mathbf{i}_{I^{c}}}
\right\|_{\mathcal{K}} .
\end{equation}
Here, when \(q>d\), the set \(I\subseteq[d]\) is identified with the same subset of
\([q]\).

It is enough to prove the localized estimate
\begin{equation}\label{A.2}
\left\|
(R^{(d)}\circ \mathfrak{p}^{(d)}\circ \mathbf 1_{L(\mathcal I)})_{\mathbf{i}_{I^{c}}}
\right\|_\mathcal{J}
\lesssim_{\alpha,D}
\sum_{q=d+1}^{D}
\sum_{\substack{\mathcal{K}\in \Pi([q]\setminus I)\\ |\mathcal{K}|\ge |\mathcal{J}|}}
\left\|
(A_q\circ \mathfrak{p}^{(q)})_{\mathbf{i}_{I^{c}}}
\right\|_{\mathcal{K}}
\end{equation}
for every \(\mathcal I\in \Pi_D\), since
\[
R^{(d)}
=
\sum_{\mathcal I\in \Pi_D}
R^{(d)}\circ \mathbf 1_{L(\mathcal I)}
\]
and the number of partitions of \([d]\) depends only on \(d\).

Fix
\[
\mathcal I=\{I_1,\ldots,I_\nu\}\in \Pi_d,
\qquad
l_\beta:=|I_\beta|,\quad \beta=1,\ldots,\nu.
\]
For \(\mathbf{i}\in L(\mathcal I)\), denote by \(j_\beta=j_\beta(\mathbf{i})\) the common value of
the coordinates $\{i_r, r\in I_\beta\}$.

Let
\[
\mathbf k=(k_1,\ldots,k_\nu)\in\mathbb N_+^\nu
\]
satisfy
\[
k_\beta\ge l_\beta\quad\text{for all }\beta\in [\nu],\qquad
\mathbf k\neq (l_1,\ldots,l_\nu),\qquad
|\mathbf k|:=k_1+\cdots+k_\nu\le D .
\]
In particular, \(q=|\mathbf k|>d\).

Define the \(d\)-tensor
\[
S_{\mathcal I}^{(d,\mathbf k)}
=
\left(s_{\mathbf{i}}^{(d,\mathbf k)}\right)_{\mathbf{i}\in[n]^d}
\]
by
\[
s_{\mathbf{i}}^{(d,\mathbf k)}
=
\mathbf 1_{\{\mathbf{i}\in L(\mathcal I)\}}\,
a^{\mathbf k}_{j_1,\ldots,j_\nu}
\prod_{\beta=1}^{\nu}
\mathbb E X_{j_\beta}^{k_\beta-l_\beta},
\]
where \(j_\beta\) is the value of \(\mathbf{i}\) on the level set \(I_\beta\), and
\(a^{\mathbf k}_{j_1,\ldots,j_\nu}\) denotes the coefficient tensor used in the
definition of \(A_q\) on level sets of sizes \(k_1,\ldots,k_\nu\).

By the definition of \(R^{(d)}\), we have
\begin{equation}\label{A.3}
R^{(d)}\circ \mathbf 1_{L(\mathcal I)}
=
\sum_{\substack{\mathbf k\in\mathbb N_+^\nu\\
k_\beta\ge l_\beta,\ \mathbf k\neq(l_1,\ldots,l_\nu),\ |\mathbf k|\le D}}
\nu!
\prod_{\beta=1}^{\nu}
\frac{k_\beta!}{(k_\beta-l_\beta)!}
S_{\mathcal I}^{(d,\mathbf k)} .
\end{equation}
The sum is finite and all combinatorial coefficients are bounded by constants
depending only on \(D\). Hence it is enough to prove, for every such
\(\mathbf k\),
\begin{equation}\label{A.4}
\left\|
(S_{\mathcal I}^{(d,\mathbf k)}\circ \mathfrak{p}^{(d)})_{\mathbf{i}_{I^{c}}}
\right\|_\mathcal{J}
\lesssim_{\alpha,D}
\sum_{\substack{\mathcal{K}\in \Pi([q]\setminus I)\\ |\mathcal{K}|\ge |\mathcal{J}|}}
\left\|
(A_q\circ \mathfrak{p}^{(q)})_{\mathbf{i}_{I^{c}}}
\right\|_{\mathcal{K}} ,
\qquad q=|\mathbf k|.
\end{equation}

Fix such a vector \(\mathbf k\). Choose a partition
\[
\widetilde{\mathcal I}
=
\{\widetilde I_1,\ldots,\widetilde I_\nu\}\in\Pi_q
\]
such that
\[
I_\beta\subseteq \widetilde I_\beta,
\qquad
|\widetilde I_\beta|=k_\beta,
\qquad
\beta=1,\ldots,\nu.
\]
For each \(\beta\), let
\[
\rho_\beta:=\min I_\beta .
\]
Define a map \(\pi:[q]\to[d]\) by
\[
\pi(r)=
\begin{cases}
r, & r\le d,\\
\rho_\beta, & r\in \widetilde I_\beta\setminus I_\beta .
\end{cases}
\]
Thus, if \(u\in L(\widetilde{\mathcal I})\), then
\[
u_r=u_{\pi(r)},\qquad r=d+1,\ldots,q.
\]

Define the \(q\)-tensor
\[
\widetilde S^{(q,\mathbf k)}
=
\left(\widetilde s_u^{(q,\mathbf k)}\right)_{u\in[n]^q}
\]
by
\begin{equation}\label{A.5}
\widetilde s_u^{(q,\mathbf k)}
=
s^{(d,\mathbf k)}_{u_{[d]}}\,
\mathbf 1_{\{u\in L(\widetilde{\mathcal I})\}} .
\end{equation}
By construction of \(A_q\), for \(u\in L(\widetilde{\mathcal I})\),
\begin{equation}\label{A.6}
\widetilde s_u^{(q,\mathbf k)}
=
(A_q)_u
\prod_{\beta=1}^{\nu}
\mathbb E X_{u_{\rho_\beta}}^{k_\beta-l_\beta}.
\end{equation}
Moreover, on \(L(\widetilde{\mathcal I})\), the coordinates \(d+1,\ldots,q\) only
repeat old coordinates. Therefore the set of distinct values of \(u\) is the same as
the set of distinct values of \(u_{[d]}\), and hence
\begin{equation}\label{A.7}
\mathfrak{p}^{(q)}_u=\mathfrak{p}^{(d)}_{u_{[d]}},
\qquad u\in L(\widetilde{\mathcal I}).
\end{equation}

We now construct a partition of \([q]\setminus I\) from \(\mathcal{J}\). Write
\[
\mathcal{J}=\{J_1,\ldots,J_\mu\}\in \Pi([d]\setminus I).
\]
Initialize
\[
K_a:=J_a,\qquad a=1,\ldots,\mu,
\qquad
K_{\mu+1}:=\varnothing .
\]
For every \(r=d+1,\ldots,q\), let \(t=\pi(r)\in[d]\). If \(t\in[d]\setminus I\),
let \(a\in\{1,\ldots,\mu\}\) be the unique index such that \(t\in J_a\), and put
\[
K_a:=K_a\cup\{r\}.
\]
If \(t\in I\), put
\[
K_{\mu+1}:=K_{\mu+1}\cup\{r\}.
\]
If \(K_{\mu+1}=\varnothing\), we discard it. The resulting family, still denoted by
$
\mathcal{K},
$
is a partition of \([q]\setminus I\), and
\begin{equation}\label{A.8}
|\mathcal{K}|\in\{|\mathcal{J}|,|\mathcal{J}|+1\}.
\end{equation}

We claim that
\begin{equation}\label{A.9}
\left\|
(S_{\mathcal I}^{(d,\mathbf k)}\circ \mathfrak{p}^{(d)})_{\mathbf{i}_{I^{c}}}
\right\|_\mathcal{J}
\le
\left\|
(\widetilde S^{(q,\mathbf k)}\circ \mathfrak{p}^{(q)})_{\mathbf{i}_{I^{c}}}
\right\|_{\mathcal{K}} .
\end{equation}

Let
\[
x^{(a)}
=
\left(x^{(a)}_{i_{J_a}}\right)_{i_{J_a}},
\qquad a=1,\ldots,\mu,
\]
be test tensors with
\[
\|x^{(a)}\|_{\mathrm{HS}}\le 1 .
\]
We define test tensors for the partition \(\mathcal{K}\). For \(a=1,\ldots,\mu\), set
\[
y^{(a)}_{u_{K_a}}
=
x^{(a)}_{u_{J_a}}
\prod_{r\in K_a\setminus[d]}
\mathbf 1_{\{u_r=u_{\pi(r)}\}} .
\]
This is well-defined because if \(r\in K_a\setminus[d]\), then
\(\pi(r)\in J_a\subseteq K_a\). Also,
\[
\|y^{(a)}\|_{\mathrm{HS}}\le \|x^{(a)}\|_{\mathrm{HS}}\le 1,
\]
since the indicator factors determine the extra coordinates \(K_a\setminus[d]\)
once \(u_{J_a}\) is fixed.

If \(K_{\mu+1}\neq\varnothing\), define
\[
y^{(\mu+1)}_{u_{K_{\mu+1}}}
=
\prod_{r\in K_{\mu+1}}
\mathbf 1_{\{u_r=i_{\pi(r)}\}} .
\]
Then
\[
\|y^{(\mu+1)}\|_{\mathrm{HS}}\le 1,
\]
because this tensor has at most one nonzero entry.

Using \eqref{A.5} and \eqref{A.7}, we obtain
\[
\begin{aligned}
&\sum_{i_{[d]\setminus I}}
(S_{\mathcal I}^{(d,\mathbf k)}\circ \mathfrak{p}^{(d)})_i
\prod_{a=1}^{\mu}x^{(a)}_{i_{J_a}}
\\
&\qquad =
\sum_{u_{[q]\setminus I}}
(\widetilde S^{(q,\mathbf k)}\circ \mathfrak{p}^{(q)})_u
\prod_{a=1}^{\mu}y^{(a)}_{u_{K_a}}
\cdot
\begin{cases}
1, & K_{\mu+1}=\varnothing,\\
y^{(\mu+1)}_{u_{K_{\mu+1}}}, & K_{\mu+1}\neq\varnothing.
\end{cases}
\end{aligned}
\]
Taking the supremum over all admissible \(x^{(a)}\)'s proves \eqref{A.9}.

It remains to compare \(\widetilde S^{(q,\mathbf k)}\circ \mathfrak{p}^{(q)}\) with
\(A_q\circ \mathfrak{p}^{(q)}\). For \(r=1,\ldots,q\), define vectors \(v_r\in\mathbb R^n\) by
\[
v_{\rho_\beta}(a)
=
\mathbb E X_a^{k_\beta-l_\beta},
\qquad a\in[n],\quad \beta=1,\ldots,\nu,
\]
and
\[
v_r(a)=1
\]
for all remaining \(r\). Since \(\|X_a\|_{\Psi_\alpha}\le 1\) and
\(0\le k_\beta-l_\beta\le D\), we have
\[
\|v_r\|_\infty\lesssim_{\alpha,D}1,
\qquad r=1,\ldots,q.
\]
By \eqref{A.6},
\[
\widetilde S^{(q,\mathbf k)}\circ \mathfrak{p}^{(q)}
=
\left(A_q\circ \mathfrak{p}^{(q)}\circ \mathbf 1_{L(\widetilde{\mathcal I})}\right)
\circ
\bigotimes_{r=1}^{q}v_r .
\]
Therefore, applying Lemma 3.4 (5) in \cite{gotze2021concentration} to the tensor
\(A_q\circ \mathfrak{p}^{(q)}\circ \mathbf 1_{L(\widetilde{\mathcal I})}\), we get
\[
\left\|
(\widetilde S^{(q,\mathbf k)}\circ \mathfrak{p}^{(q)})_{\mathbf{i}_{I^{c}}}
\right\|_{\mathcal{K}}
\lesssim_{\alpha,D}
\left\|
(A_q\circ \mathfrak{p}^{(q)}\circ \mathbf 1_{L(\widetilde{\mathcal I})})_{\mathbf{i}_{I^{c}}}
\right\|_{\mathcal{K}} .
\]
Then Lemma \ref{Lem_matrixnorm_comparasion} gives
\[
\left\|
(A_q\circ \mathfrak{p}^{(q)}\circ \mathbf 1_{L(\widetilde{\mathcal I})})_{\mathbf{i}_{I^{c}}}
\right\|_{\mathcal{K}}
\lesssim_D
\left\|
(A_q\circ \mathfrak{p}^{(q)})_{\mathbf{i}_{I^{c}}}
\right\|_{\mathcal{K}} .
\]
Combining the last two displays with \eqref{A.9}, we obtain
\[
\left\|
(S_{\mathcal I}^{(d,\mathbf k)}\circ \mathfrak{p}^{(d)})_{\mathbf{i}_{I^{c}}}
\right\|_\mathcal{J}
\lesssim_{\alpha,D}
\left\|
(A_q\circ \mathfrak{p}^{(q)})_{\mathbf{i}_{I^{c}}}
\right\|_{\mathcal{K}} .
\]
By \eqref{A.8}, this implies
\[
\left\|
(S_{\mathcal I}^{(d,\mathbf k)}\circ \mathfrak{p}^{(d)})_{\mathbf{i}_{I^{c}}}
\right\|_\mathcal{J}
\lesssim_{\alpha,D}
\sum_{\substack{\mathcal{K}'\in \Pi([q]\setminus I)\\ |\mathcal{K}'|\ge |\mathcal{J}|}}
\left\|
(A_q\circ \mathfrak{p}^{(q)})_{\mathbf{i}_{I^{c}}}
\right\|_{\mathcal{K}'} .
\]
This proves \eqref{A.4}.

Finally, summing over the finitely many admissible \(\mathbf k\)'s in
\eqref{A.3}, and then over the finitely many level-set partitions
\(\mathcal I\in \Pi_D\), proves \eqref{A.1}. Hence \((4.5)\) follows.

\section{Supplementary estimates for Corollary~\ref{cor:cor1_corrected_refined}}
\label{app:cor1_appendixB}

In this appendix we collect the deterministic tensor-norm estimates used in the
proof of Corollary~\ref{cor:cor1_corrected_refined}.  
Throughout this appendix, set
\[
    \mathsf L:=[d],
    \qquad
    \mathsf R:=d+[d]=\{d+1,\ldots,2d\}.
\]
Thus
\[
    \|A\|_{\mathsf L,\mathsf R}
    =
    \|A\|_{[d],[2d]\setminus[d]} .
\]

\subsection*{The unweighted trace tensors}

Let
\[
    I=\{k_1,\ldots,k_l\}\subseteq[d],
    \qquad k_1<\cdots<k_l,
\]
and let $a=(a_1,\ldots,a_l), b=(b_1,\ldots,b_l)\in[n]^l$.  For
$u=(u_t)_{t\in[d]\setminus I}\in[n]^{d-l}$ and
$\varepsilon=(\varepsilon_1,\ldots,\varepsilon_l)\in\{0,1\}^l$, recall the
filling map
\[
    \Theta_{I,\varepsilon}(a,b,u)=(i_1,\ldots,i_{2d})\in[n]^{2d},
\]
defined by
\[
    i_t=i_{d+t}=u_t,
    \qquad t\in[d]\setminus I,
\]
and, for $t=k_m\in I$,
\[
    (i_{k_m},i_{d+k_m})
    =
    \begin{cases}
        (a_m,b_m), & \varepsilon_m=0,\\
        (b_m,a_m), & \varepsilon_m=1.
    \end{cases}
\]
We define the unweighted partial trace associated with one selected-pair swap
pattern by
\begin{equation}\label{eq:B_T_I_epsilon_A}
    T_{I,\varepsilon}A(a,b)
    :=
    \sum_{u\in[n]^{d-l}} A_{\Theta_{I,\varepsilon}(a,b,u)},
    \qquad a,b\in[n]^l.
\end{equation}
In the homogeneous case $p_1=\cdots=p_n=p$, the  reduced tensor from
Theorem~\ref{thm:tensor_poly_concentration} can be written as
\begin{equation}\label{eq:B_corrected_tensor_homogeneous_decomposition}
    \mathcal A_I^{(2l,p)}
    =
    p^{d-l}\,\widetilde{\mathfrak{P}}^{(2l)}
    \circ
    \sum_{\varepsilon\in\{0,1\}^l}T_{I,\varepsilon}A.
\end{equation}

For later use, define the bijection
\[
    \rho_{I,\varepsilon}:[2l]\to I\cup(I+d)
\]
by
\begin{equation}\label{eq:B_rho_I_epsilon}
    \rho_{I,\varepsilon}(m)
    =
    \begin{cases}
        k_m, & \varepsilon_m=0,\\
        d+k_m, & \varepsilon_m=1,
    \end{cases}
    \qquad
    \rho_{I,\varepsilon}(l+m)
    =
    \begin{cases}
        d+k_m, & \varepsilon_m=0,\\
        k_m, & \varepsilon_m=1,
    \end{cases}
    \quad m\in[l].
\end{equation}
If $\mathcal J$ is a partition of a subset of $[2l]$, then
$\rho_{I,\varepsilon}(\mathcal J)$ denotes the partition obtained by applying
$\rho_{I,\varepsilon}$ to every block of $\mathcal J$.

\begin{lemma}[Lemma 23 in \cite{BKW22}]
\label{lem:B_block_splitting}
Let $B=(B_i)_{i\in[n]^m}$ be an $m$-tensor and let
$\mathcal K=\{K_1,\ldots,K_s\}$ be a partition of $[m]$. Suppose that
$K_s=K_s'\cup K_{s+1}'$ is a disjoint union. Then
\begin{equation}\label{eq:B_block_splitting}
    \|B\|_{K_1,\ldots,K_{s-1},K_s',K_{s+1}'}
    \le
    \|B\|_{K_1,\ldots,K_s}
    \le
    n^{\min\{|K_s'|,|K_{s+1}'|\}/2}
    \|B\|_{K_1,\ldots,K_{s-1},K_s',K_{s+1}'}.
\end{equation}
\end{lemma}

\begin{lemma}
\label{lem:B_free_weight_multiplier}
Assume $p_1=\cdots=p_n=p\in(0,1]$.  Let $l\in[d]$ and let
$B=(B_j)_{j\in[n]^{2l}}$ be a $2l$-tensor.

\begin{enumerate}
    \item For every partition $\mathcal J\in\Pi([2l])$,
    \begin{equation}\label{eq:B_free_weight_multiplier_empty}
        \bigl\|\widetilde{\mathfrak{P}}^{(2l)}\circ B\bigr\|_{\mathcal J}
        \lesssim_d
        p^{l/2}\|B\|_{\mathcal J}.
    \end{equation}

    \item Let $I_2\subseteq[2l]$, fix $j_{I_2}\in[n]^{I_2}$, and let
    $\mathcal J\in\Pi([2l]\setminus I_2)$. Then
    \begin{equation}\label{eq:B_free_weight_multiplier_sliced}
        \frac{
        \bigl\|
        (\widetilde{\mathfrak{P}}^{(2l)}\circ B)_{j_{I_{2}^{c}}}
        \bigr\|_{\mathcal J}
        }{
        \widetilde{\mathfrak{P}}_{I_2}^{(2l)}(j_{I_2})
        }
        \lesssim_d
        \bigl\|B_{j_{I_{2}^{c}}}\bigr\|_{\mathcal J}.
    \end{equation}
\end{enumerate}
\end{lemma}

\begin{proof}
Write the multi-index as $j=(a,b)=(a_1,\ldots,a_l,b_1,\ldots,b_l)$. In the homogeneous case, the global weight tensor factors completely over the block structure:
\[
    \widetilde{\mathfrak{P}}^{(2l)}(a,b) = \prod_{m=1}^l q(a_m,b_m),
\]
where the pairwise weight function is given by
\[
    q(a_m,b_m) = p\,\mathbf 1_{\{a_m\ne b_m\}} + \sqrt p\,\mathbf 1_{\{a_m=b_m\}} = p+(\sqrt p-p)\mathbf 1_{\{a_m=b_m\}}.
\]
Expanding the product over $m=1,\dots,l$, the weight tensor can be expressed as a finite linear combination of indicator functions associated with specific equality patterns:
\[
    \widetilde{\mathfrak{P}}^{(2l)}(a,b) = \sum_{S\subseteq[l]} p^{l-|S|}(\sqrt p-p)^{|S|} \mathbf 1_{E_S}(a,b),
\]
where $E_S = \{ (a,b) : a_m = b_m \text{ for all } m \in S \}$. Geometrically, each set $E_S$ is an intersection of generalized diagonals. According to Lemma 3.4 in \cite{gotze2021concentration}, the operation of taking the Hadamard product with the indicator function of a generalized diagonal (or intersections thereof) constitutes a bounded multiplier on any partition norm $\|\cdot\|_{\mathcal J}$. Therefore, there exists a constant $C_d > 0$ depending only on the tensor dimension $d$ such that $\|\mathbf 1_{E_S} \circ B\|_{\mathcal J} \le C_d \|B\|_{\mathcal J}$ for all $S \subseteq [l]$.

Since $p \in (0,1]$, we have $\sqrt p - p \ge 0$. The total variation (the sum of the absolute values of the coefficients) converges exactly to $p^{l/2}$ by the binomial theorem:
\[
    \sum_{S\subseteq[l]} \left| p^{l-|S|}(\sqrt p-p)^{|S|} \right| = \sum_{S\subseteq[l]} p^{l-|S|}(\sqrt p-p)^{|S|} = \bigl(p + (\sqrt p-p)\bigr)^l = p^{l/2}.
\]
Applying the triangle inequality for the partition norm gives:
\[
    \bigl\|\widetilde{\mathfrak{P}}^{(2l)}\circ B\bigr\|_{\mathcal J} \le \sum_{S\subseteq[l]} p^{l-|S|}(\sqrt p-p)^{|S|} \bigl\|\mathbf 1_{E_S} \circ B\bigr\|_{\mathcal J}  \lesssim_d p^{l/2}\|B\|_{\mathcal J}.
\]
This establishes \eqref{eq:B_free_weight_multiplier_empty}.

For the sliced estimate \eqref{eq:B_free_weight_multiplier_sliced}, we fix a subset of coordinates $I_2 \subseteq [2l]$ with values $j_{I_2}$. The condition weight multiplier applied to the remaining free coordinates $j_{I_2^{\,c}}$ is the quotient
\[
    W(j_{I_2^{\,c}}) = \frac{\widetilde{\mathfrak{P}}^{(2l)}(j)}{\widetilde{\mathfrak{P}}_{I_2}^{(2l)}(j_{I_2})}.
\]
Due to the independence of the blocks, this quotient preserves the product structure, yielding $W(j_{I_2^{\,c}}) = \prod_{m=1}^l w_m$. For any given block $m$, the factor $w_m$ falls into one of three mutually exclusive cases depending on the fixed coordinates:
\begin{enumerate}
    \item \textbf{Neither coordinate is fixed} ($a_m, b_m \notin I_2$): $w_m = q(a_m, b_m) = p + (\sqrt p - p)\mathbf 1_{\{a_m=b_m\}}$. The sum of absolute coefficients is $p + (\sqrt p - p) = \sqrt p \le 1$.
    \item \textbf{Exactly one coordinate is fixed} (e.g., $b_m = x_0 \in j_{I_2}$ is fixed and $a_m = x$ is free): Dividing the joint weight by the fixed marginal weight yields $w_m = \sqrt p + (1-\sqrt p)\mathbf 1_{\{x=x_0\}}$. The sum of absolute coefficients is $\sqrt p + (1-\sqrt p) = 1$.
    \item \textbf{Both coordinates are fixed} ($a_m, b_m \in I_2$): The factor identically cancels out, yielding $w_m = 1$. The sum of absolute coefficients is $1$.
\end{enumerate}
Expanding the product $W(j_{I_2^{\,c}}) = \prod_{m=1}^l w_m$ generates a new linear combination $W(j_{I_2^{\,c}}) = \sum_K \tilde{c}_K \mathbf 1_{E'_K}$, where each $E'_K$ represents either an equality pattern strictly among the free coordinates, or an equality forcing a free coordinate to match a fixed value (which corresponds to a coordinate projection after slicing). 
The total variation of this expansion is simply the product of the absolute coefficient sums of the individual blocks:
\[
    \sum_K |\tilde{c}_K| \lesssim_d 1.
\]
Applying the triangle inequality and Lemma 3.4 in \cite{gotze2021concentration} yields:
\[
    \frac{ \bigl\| (\widetilde{\mathfrak{P}}^{(2l)}\circ B)_{j_{I_2^{\,c}}} \bigr\|_{\mathcal J} }{ \widetilde{\mathfrak{P}}_{I_2}^{(2l)}(j_{I_2}) } = \biggl\| \sum_K \tilde{c}_K \bigl( \mathbf 1_{E'_K} \circ B_{j_{I_2^{\,c}}} \bigr) \biggr\|_{\mathcal J}  \lesssim_d \bigl\|B_{j_{I_2^{\,c}}}\bigr\|_{\mathcal J}.
\]
This establishes \eqref{eq:B_free_weight_multiplier_sliced} and completes the proof.
\end{proof}

\begin{lemma}[Sliced partition norms versus the middle flattening]
\label{lem:B_sliced_partition_vs_middle}
Let $F\subseteq[2d]$ with $|F|=q$, and fix $j_F\in[n]^F$. Let
$\mathcal K=\{K_1,\ldots,K_m\}$ be a partition of
$F^c=[2d]\setminus F$. Then
\begin{equation}\label{eq:B_sliced_partition_vs_middle}
    \|A_{j_F^{\,c}}\|_{\mathcal K}
    \le
    n^{d-(m+q)/2}\|A\|_{\mathsf L,\mathsf R}.
\end{equation}
\end{lemma}

\begin{proof}
For each block $K_s\in\mathcal K$, set
\[
    K_s^-:=K_s\cap\mathsf L,
    \qquad
    K_s^+:=K_s\cap\mathsf R.
\]
Let $\mathcal K'$ be the partition obtained from $\mathcal K$ by replacing
every nonempty mixed block $K_s$ by the two nonempty blocks $K_s^-$ and
$K_s^+$. Repeatedly applying Lemma~\ref{lem:B_block_splitting} yields
\begin{equation}\label{eq:B_split_all_mixed_blocks}
    \|A_{j_F^{\,c}}\|_{\mathcal K}
    \le
    n^{\frac12\sum_{s=1}^m\min\{|K_s^-|,|K_s^+|\}}
    \|A_{j_F^{\,c}}\|_{\mathcal K'}.
\end{equation}
The partition $\mathcal K'$ refines the two-block partition
\[
    \{\mathsf L\setminus F,\,\mathsf R\setminus F\},
\]
with empty blocks omitted. By monotonicity of partition norms under refinement,
and by embedding the test vectors into the fixed coordinates $j_F$, we have
\begin{equation}\label{eq:B_refinement_to_middle}
    \|A_{j_F^{\,c}}\|_{\mathcal K'}
    \le
    \|A_{j_F^{\,c}}\|_{\mathsf L\setminus F,\,\mathsf R\setminus F}
    \le
    \|A\|_{\mathsf L,\mathsf R}.
\end{equation}
Finally,
\[
    \sum_{s=1}^m\min\{|K_s^-|,|K_s^+|\}
    \le
    \sum_{s=1}^m(|K_s|-1)
    =
    |F^c|-m
    =
    2d-q-m.
\]
Combining this estimate with \eqref{eq:B_split_all_mixed_blocks} and
\eqref{eq:B_refinement_to_middle} proves \eqref{eq:B_sliced_partition_vs_middle}.
\end{proof}

\begin{lemma}
\label{lem:B_trace_lifting}
\label{Lem_norm_estimation}
Let $A=(A_{i_1,\ldots,i_{2d}})$ be a $2d$-tensor. Fix
$I=\{k_1,\ldots,k_l\}\subseteq[d]$ with $k_1<\cdots<k_l$, and fix
$\varepsilon\in\{0,1\}^l$. Let $T_{I,\varepsilon}A$ be defined by
\eqref{eq:B_T_I_epsilon_A}, and let $\rho_{I,\varepsilon}$ be defined by
\eqref{eq:B_rho_I_epsilon}.

Let $I_2\subseteq[2l]$ and set
\[
    F:=\rho_{I,\varepsilon}(I_2)\subseteq[2d].
\]
For a fixed $j_{I_2}$, denote by $j_F^\rho$ the corresponding fixed coordinates
of $A$, namely
\[
    i_{\rho_{I,\varepsilon}(s)}=j_s,
    \qquad s\in I_2.
\]
If $\mathcal J=\{J_1,\ldots,J_\kappa\}\in
\Pi([2l]\setminus I_2)$, define
\begin{equation}\label{eq:B_K_from_J_rho}
    \mathcal K
    :=
    \rho_{I,\varepsilon}(\mathcal J)
    \cup
    \bigl\{\{t,d+t\}:t\in[d]\setminus I\bigr\}.
\end{equation}
Then $\mathcal K$ is a partition of $[2d]\setminus F$, and
\begin{equation}\label{eq:B_trace_lifting_sliced}
    \bigl\|
    (T_{I,\varepsilon}A)_{j_{I_{2}^{c}}}
    \bigr\|_{\mathcal J}
    \le
    n^{(d-l)/2}
    \bigl\|A_{(j_F^\rho)^{\,c}}\bigr\|_{\mathcal K}.
\end{equation}
In particular, if $I_2=\varnothing$, then
\begin{equation}\label{eq:B_trace_lifting_unsliced}
    \|T_{I,\varepsilon}A\|_{\mathcal J}
    \le
    n^{(d-l)/2}
    \|A\|_{\rho_{I,\varepsilon}(\mathcal J),\,\{\{t,d+t\}:t\in[d]\setminus I\}}.
\end{equation}
\end{lemma}

\begin{proof}
It suffices to prove \eqref{eq:B_trace_lifting_sliced}. Let
$x^{(s)}=(x^{(s)}_{j_{J_s}})$, $s=1,\ldots,\kappa$, be test tensors with
Hilbert--Schmidt norm at most one. By the definition of $T_{I,\varepsilon}A$,
\begin{align*}
&\sum_{j_{[2l]\setminus I_2}}
(T_{I,\varepsilon}A)(j)
\prod_{s=1}^{\kappa}x^{(s)}_{j_{J_s}}                                               \\
&\qquad=
\sum_{i_{[2d]\setminus F}}
A_i
\prod_{s=1}^{\kappa}x^{(s)}_{i_{\rho_{I,\varepsilon}(J_s)}}
\prod_{t\in[d]\setminus I}\mathbf 1_{\{i_t=i_{d+t}\}},
\end{align*}
where the coordinates in $F$ are fixed according to $j_F^\rho$.  For every
$t\in[d]\setminus I$, define
\[
    y^{(t)}_{a,b}:=n^{-1/2}\mathbf 1_{\{a=b\}},
    \qquad a,b\in[n].
\]
Then $\|y^{(t)}\|_2=1$. The last display equals
\[
    n^{(d-l)/2}
    \sum_{i_{[2d]\setminus F}}
    A_i
    \prod_{s=1}^{\kappa}x^{(s)}_{i_{\rho_{I,\varepsilon}(J_s)}}
    \prod_{t\in[d]\setminus I}y^{(t)}_{i_t,i_{d+t}}.
\]
This is bounded by
$n^{(d-l)/2}\|A_{(j_F^\rho)^c}\|_{\mathcal K}$ by the definition of the
partition norm. Taking the supremum over all admissible test tensors
$x^{(s)}$ proves \eqref{eq:B_trace_lifting_sliced}.
\end{proof}

\begin{lemma}[Uniform bounds for the unweighted partial traces]
\label{lem:B_partial_trace_estimates}
Let $I\subseteq[d]$, $|I|=l$, let $\varepsilon\in\{0,1\}^l$, and let
$T_{I,\varepsilon}A$ be defined by \eqref{eq:B_T_I_epsilon_A}. Let
$I_2\subseteq[2l]$, set $q:=|I_2|$, and let
$\mathcal J\in\Pi([2l]\setminus I_2)$ with
$\kappa:=|\mathcal J|$. Then, uniformly in the fixed coordinates $j_{I_2}$,
\begin{equation}\label{eq:B_partial_trace_op_estimate}
    \bigl\|
    (T_{I,\varepsilon}A)_{j_{I_{2}^{c}}}
    \bigr\|_{\mathcal J}
    \le
    n^{d-(\kappa+q)/2}\|A\|_{\mathsf L,\mathsf R}.
\end{equation}
Moreover, when $I_2=\varnothing$,
\begin{equation}\label{eq:B_partial_trace_hs_estimate}
    \|T_{I,\varepsilon}A\|_{\mathcal J}
    \le
    n^{(d-l)/2}\|A\|_{[2d]}.
\end{equation}
\end{lemma}

\begin{proof}
Let $\rho=\rho_{I,\varepsilon}$ and set $F:=\rho(I_2)$. Let
\[
    \mathcal K
    :=
    \rho(\mathcal J)
    \cup
    \bigl\{\{t,d+t\}:t\in[d]\setminus I\bigr\}.
\]
Then $\mathcal K$ is a partition of $F^c$ and
\[
    |\mathcal K|=\kappa+d-l.
\]
By Lemma~\ref{lem:B_trace_lifting},
\begin{equation}\label{eq:B_trace_to_full_tensor_norm_new}
    \bigl\|
    (T_{I,\varepsilon}A)_{j_{I_{2}^{c}}}
    \bigr\|_{\mathcal J}
    \le
    n^{(d-l)/2}\|A_{(j_F^\rho)^{\,c}}\|_{\mathcal K}.
\end{equation}
Applying Lemma~\ref{lem:B_sliced_partition_vs_middle} with
$m=|\mathcal K|=\kappa+d-l$ and $q=|F|=|I_2|$ gives
\[
    \|A_{(j_F^\rho)^{\,c}}\|_{\mathcal K}
    \le
    n^{d-(\kappa+d-l+q)/2}\|A\|_{\mathsf L,\mathsf R}.
\]
Multiplying by $n^{(d-l)/2}$ proves \eqref{eq:B_partial_trace_op_estimate}.

For \eqref{eq:B_partial_trace_hs_estimate}, take $I_2=\varnothing$ in
\eqref{eq:B_trace_to_full_tensor_norm_new}. Then $F=\varnothing$ and
$\mathcal K$ is a partition of $[2d]$. Since every partition norm is bounded by
the Hilbert--Schmidt norm,
\[
    \|A\|_{\mathcal K}\le \|A\|_{[2d]},
\]
and \eqref{eq:B_partial_trace_hs_estimate} follows.
\end{proof}

\paragraph{Acknowledgment}
K. Wang was partially supported by the Hong Kong RGC grant GRF 16304222. The authors would like to thank Rados\l{}aw Adamczak, Yiyun He and Yizhe Zhu for inspiring and valuable discussions.

\bibliography{Polynomial}
\bibliographystyle{abbrv}

\end{document}